\documentclass[11pt]{article}
\usepackage{amssymb,amsmath,latexsym, epsfig}

\begin{document}

\begin{doublespace}

\newtheorem{thm}{Theorem}[section]
\newtheorem{lemma}[thm]{Lemma}
\newtheorem{defn}{Definition}[section]
\newtheorem{prop}[thm]{Proposition}
\newtheorem{corollary}[thm]{Corollary}
\newtheorem{remark}[thm]{Remark}
\newtheorem{example}[thm]{Example}
\numberwithin{equation}{section}

\def\ee{\varepsilon}
\def\qed{{\hfill $\Box$ \bigskip}}
\def\MM{{\cal M}}
\def\BB{{\cal B}}
\def\LL{{\cal L}}
\def\FF{{\cal F}}
\def\GG{{\cal G}}
\def\EE{{\cal E}}
\def\QQ{{\cal Q}}

\def\R{{\mathbb R}}
\def\L{{\bf L}}
\def\E{{\mathbb E}}
\def\F{{\bf F}}
\def\P{{\mathbb P}}
\def\N{{\mathbb N}}
\def\eps{\varepsilon}
\def\wh{\widehat}
\def\pf{\noindent{\bf Proof.} }

\title{\Large \bf On the Estimates of the Density of Feynman-Kac Semigroups of $\alpha$-Stable-like Processes }
\author{ Chunlin Wang\\
Department of Mathematics\\
University of Illinois\\
Urbana, IL 61801, USA\\
Email: cwang13@uiuc.edu}
\date{}
\maketitle

\begin{abstract}
 Suppose that $\alpha \in (0,2)$ and that $X$ is an
$\alpha$-stable-like process on $\R^d$. Let $F$ be a function on
$\R^d$ belonging to the class $\bf{J_{d,\alpha}}$ (see
Introduction) and $A_{t}^{F}$ be $\sum_{s \le
t}F(X_{s-},X_{s}),\,\, t> 0$, a discontinuous additive functional
of $X$. With neither $F$ nor $X$ being symmetric, under certain
conditions, we show that the Feynman-Kac semigroup $\{S_{t}^{F}:t
\ge 0\}$ defined by
$$ S_{t}^{F}f(x)=\mathbb{E}_{x}(e^{-A_{t}^{F}}f(X_{t}))$$
has a density $q$ and that there exist positive constants
$C_1,C_2,C_3$ and $C_4$ such that
$$C_{1}e^{-C_{2}t}t^{-\frac{d}{\alpha}}\left(1 \wedge \frac{t^{\frac{1}{\alpha}}}{|x-y|}\right)^{d+\alpha} \leq q(t,x,y) \leq  C_{3}e^{C_{4}t}t^{-\frac{d}{\alpha}}\left(1 \wedge \frac{t^{\frac{1}{\alpha}}}{|x-y|}\right)^{d+\alpha}$$
for all $(t,x,y)\in (0,\infty) \times \R^d \times \R^d$.

\end{abstract}

\noindent {\bf AMS 2000 Mathematics Subject Classification}:
Primary 60J45, 60J35; Secondary 60J55, 60J75, 60G51, 60G52

\noindent{\bf Keywords and phrases:} $\alpha$-stable-like
processes, nonsymmetric, Kato class, Feynman-Kac semigroups,
discontinuous additive functionals

 \vspace{.1truein}
\noindent{\bf Running Title:} Density estimates of Feynman-Kac
semigroups of stable processes

\pagebreak
\section{Introduction}
Suppose $X=(X_t,\mathbb{P}_t)$ is a Hunt process on $\R^d$ with a
L\'{e}vy system $(N,H)$ given by $H_{t}=t$ and
$$ N(x,dy)=2C(x,y){|x-y|}^{-(d+\alpha)}m(dy),$$
where $m$ is a measure on $\R^d$ given by $m(dx)=M(x)dx$ with
$M(x)$ bounded between two positive numbers. That is for any
nonnegative function $f$ on $\R^d\times\R^d$ vanishing on the
diagonal
$$\mathbb{E}_x\left(\sum_{s \le
t}f(X_{s-},X_{s})\right)=\mathbb{E}_x\int_{0}^{t}\int_{R^d}\frac{2C(X_s,y)f(X_s,y)}{{|X_s-y|}^{d+\alpha}}m(dy)ds,$$
for every $x \in \R^d$ and $t>0$.

We introduce $\alpha$-stable-like processes as follows.
\begin{defn}
We say that $X$ is an $\alpha$-stable-like process if $C(x,y)$ is
bounded.
\end{defn}

In this paper, we assume that $X$ admits a transition density
$p(t,x,y)$ with respect to $m$ and $p(t,x,y)$ is jointly
continuously on $(0,\infty)\times\R^d\times\R^d$ and satisfies the
condition
\begin{equation}
 M_{1}t^{-\frac{d}{\alpha}}\left(1 \wedge
\frac{t^{\frac{1}{\alpha}}}{|x-y|}\right)^{d+\alpha}\leq
p(t,x,y)\leq M_{2}t^{-\frac{d}{\alpha}}\left(1 \wedge
\frac{t^{\frac{1}{\alpha}}}{|x-y|}\right)^{d+\alpha},\,\,\forall
(t,x,z)\in (0,\infty)\times\R^d\times\R^d,
\end{equation}
where $M_{1}$ and $M_{2}$ are positive constants.

Here we do not assume that $X$ is symmetric. When $X$ is
symmetric, it is called a symmetric $\alpha$-stable-like process,
which was introduced in \cite{CT}, where a symmetric Hunt process
is associated with a regular Dirichlet form and thus Dirichlet
form method can be applied. It was also shown in \cite{CT} that
the transition densities of symmetric $\alpha$-stable-like
processes satisfy (1.1).

We list some examples which are $\alpha$-stable-like processes and
satisfy (1.1). For one dimensional strictly $\alpha$-stable
processes with L\'{e}vy measure $\nu$ concentrated neither on
$(0,\infty)$ nor on $(-\infty,0)$, the L\'{e}vy measure $\nu(dx)
=c_1x^{-1-\alpha}dx$ on $(0,\infty)$ and $\nu(dx)
=c_2x^{-1-\alpha}dx$ on $(-\infty,0)$ with $c_1>0$ and $c_2>0$,
which implies $C(x,y)$ in the L\'{e}vy system as above is bounded
between two positive numbers. We set $c=c_1+c_2$ and
$\beta=(c_1-c_2)/(c_1+c_2)$. Let $\rho=(1+\beta)/2$ or
$=(1-\beta\frac{2-\alpha}{2})/2$, according to $\alpha<1$ or $>1$.
Without loss of generality, we can fix the parameter $c$ and
assume that $c$ equals $\cos(\frac{\pi\beta\alpha}{2})$,
$\frac{\pi}{2}$ or $\cos(\pi\beta\frac{2-\alpha}{2})$ for
$\alpha<1$, $=1$, or $>1$, respectively. \cite{SK} gave the
following estimates for the continuous transition density
$p(t,0,x)$, which equals $t^{-1/\alpha}p(1,0,t^{-1/\alpha}x)$:
 \begin{enumerate}
\item When $x\rightarrow \infty$,
\begin{eqnarray*}
   p(1,0,x)&\sim&\frac{1}{\pi}\Gamma(\alpha+1)(\sin(\pi
    \rho  \alpha))x^{-\alpha-1},\, \,\,\textrm{if } \alpha \neq 1,\\
   p(1,0,x)&\sim& \frac{1+\beta}{2}x^{-2}, \,\,\,\textrm{if } \alpha =1,
\end{eqnarray*}

\item When $x\rightarrow 0$,
\begin{eqnarray*}
      p(1,0,x)&\rightarrow&
      \frac{1}{\pi}\Gamma(1/\alpha+1)(\sin\pi\alpha), \,\,\,\textrm{if
      }\alpha \neq 1,\\
      p(1,0,x)&\rightarrow& \frac{1}{\pi}b_1, \,\,\,\textrm{if }\alpha=1,
      \,\beta>0,
\end{eqnarray*}
where $b_1$ is a positive constant.

 See (14.37), (14.30), (14.33)
and (14.32) in \cite{SK} for details. It is clear that the dual
process of the one dimensional strictly $\alpha$-stable process
has the transition density $p(t,0,-x)$. Thus applying the above
estimates to $p(t,0,-x)$, we get

\item When $x\rightarrow -\infty$,
\begin{eqnarray*}
   p(1,0,x)&\sim&\frac{1}{\pi}\Gamma(\alpha+1)(\sin(\pi
   \rho\alpha))|x|^{-\alpha-1},\,\,\,\textrm{if } \alpha \neq 1,\\
   p(1,0,x)&\sim& \frac{1+\beta}{2}|x|^{-2}, \,\,\,\textrm{if }
   \alpha=1.
\end{eqnarray*}
\item When $x\rightarrow 0$,
$$p(1,0,x)\rightarrow \frac{1}{\pi}\tilde{b}_1,\,\,\,\textrm{if }
\alpha=1,\, \beta<0,$$ where $\tilde{b}_1$ is a positive constant.
\end{enumerate}

One dimensional strictly $\alpha$-stable process with $\alpha=1$
and $\beta=0$ is a Cauchy process with drift $0$. It is easy to
see that when $x\rightarrow 0$, $p(1,0,x)\rightarrow$ a positive
constant.

\cite{SK} also pointed out that $p(t,0,x)$ is positive when the
L\'{e}vy measure $\nu$ is concentrated neither on $(0,\infty)$ nor
on $(-\infty,0)$.

Therefore when the L\'{e}vy measure $\nu$ is concentrated neither
on $(0,\infty)$ nor on $(-\infty,0)$, the transition density $p$
satisfies (1.1).

 For higher dimensions, \cite{VZ}
considered a large class of nonsymmetric strictly $\alpha$-stable
processes with $0<\alpha<2$, which has L\'{e}vy measure $\nu$
satisfying
$$
\nu(B)=\int_S\lambda(d\xi)\,\int^{\infty}_{0}1_{B}(r\xi)r^{-(1+\alpha)}dr$$
for every Borel set $B$ in $\R^d$, where $\lambda$ is a finite
measure on the unit sphere $S=\{x\in \R^d:\,|x|=1\}$ and is called
the spherical part of the L\'{e}vy measure $\nu$. $\lambda$ is
assumed to have a density $\phi:S\rightarrow(0,\infty)$ such that
$$\phi=\frac{d\lambda}{d\sigma} \textrm{ and } \kappa \leq
\phi(\xi)\leq \kappa^{-1},\,\,\,\forall \xi\in S,$$ where $\sigma$
is the surface measure on the unit sphere and $\kappa>0$ is a
positive constant. The assumption on $\phi$ implies the transition
density $p(t,0,x)>0$ for all $t>0$ and all $x\in \R^d$. It is
known that $p(1,0,x)$ is uniformly bounded in $x\in \R^d$.
\cite{VZ} pointed out that the L\'{e}vy measure $\nu$ has a
density $f(x)=\phi(x/|x|)|x|^{-(d+\alpha)}$ with respect to the
$d$-dimensional Lebesgue measure, and
$$\kappa|x|^{-(d+\alpha)}\leq f(x)\leq\kappa^{-1}|x|^{-(d+\alpha)}$$
for every $x\in \R^d \setminus \{0\}$. Then the transition density
$p(t,0,x)$ of the processes satisfy
$$
p(t,0,x)\leq \tilde{C}t^{-\frac{d}{\alpha}},\,
\,\,x\in\R^d,\,\,t>0,
$$
 and
$$
 p(t,0,x)\leq
\tilde{C}t|x|^{-(\alpha+d)},\,\,\,x\in\R^d\setminus \{0\},\,\,t>0,
$$
where $\tilde{C}$ is a positive constant. See (2.6) and (2.7) in
\cite{VZ} for these two inequalities. Thus we have
\begin{equation}
p(t,0,x)\leq \tilde{C}t^{-\frac{d}{\alpha}}\left(1 \wedge
\frac{t^{\frac{1}{\alpha}}}{|x-y|}\right)^{d+\alpha}      ,\,\,\,
x\in\R^d,\,\,t>0.
\end{equation}
(2.12) in \cite{VZ} gave the following estimate,
$$\tilde{c}|x|^{-(\alpha+d)}\leq p(1,0,x)\leq
\tilde{C}|x|^{-(\alpha+d)},\,\,\,\textrm{for large }x,$$ where
$\tilde{c}$ is a positive constant and $\tilde{C}$ is the same
constant as above. This implies that
$$\tilde{c}|t^{-\frac{1}{\alpha}}x|^{-(\alpha+d)}\leq p(1,0,t^{-\frac{1}{\alpha}}x)\leq
\tilde{C}|t^{-\frac{1}{\alpha}}x|^{-(\alpha+d)},\,\,\,\textrm{for
large }t^{-\frac{1}{\alpha}}x.$$ Thus
\begin{equation}
\tilde{c}t^{-\frac{d}{\alpha}}|t^{-\frac{1}{\alpha}}x|^{-(\alpha+d)}\leq
t^{-\frac{d}{\alpha}}p(1,0,t^{-\frac{1}{\alpha}}x)=p(t,0,x),\,\,\,\textrm{for
large }t^{-\frac{1}{\alpha}}x.
\end{equation}

For small $t^{-\frac{1}{\alpha}}x$, since $p(1,0,x)$ is positive
and continuous in $x\in \R^d$, there exists a positive constant
$\tilde{c}_0$ such that
$$\tilde{c}_0\leq p(1,0,t^{-\frac{1}{\alpha}}x),$$
which implies
\begin{equation}
\tilde{c}_0t^{-\frac{d}{\alpha}}\leq
t^{-\frac{d}{\alpha}}p(1,0,t^{-\frac{1}{\alpha}}x)=p(t,0,x),\,\,\,\textrm{for
small }t^{-\frac{1}{\alpha}}x.
\end{equation}

Combining (1.2), (1.3) and (1.4), we can see that the transition
density $p$ satisfies (1.1). It is clear that
$C(x,y)=\phi(\frac{x-y}{|x-y|})|x-y|^{-(d+\alpha)}$ is bounded
between two positive numbers.

 We say that a function $V$ on $\R^d$ belongs to the
Kato class $\bf{K_{d,\alpha}}$ if $$\lim_{t\downarrow 0}\sup_{x\in
\R^d}\int_{0}^{t}\int_{\R^d}p(s,x,y)|V(y)|dyds=0,$$ and we say
that a signed measure $\mu$ on $\R^d$ belongs to the Kato class
$\bf{K_{d,\alpha}}$ if $$\lim_{t\downarrow 0}\sup_{x\in
\R^d}\int_{0}^{t}\int_{\R^d}p(s,x,y)|\mu|(dy)ds=0.$$

Suppose $F$ is a function on $\R^d \times \R^d$.
\begin{defn} We
say $F$ belongs to $\bf{J_{d,\alpha}}$ if $F$ is bounded,
vanishing on the diagonal, and the function
$$
               x\mapsto \int_{\R^d}\frac{|F(x,y)|}{|x-y|^{d+\alpha}}\,dy
$$
belongs to $\bf{K_{d,\alpha}}$.
\end{defn}

For any $F \in \bf{J_{d,\alpha}}$, we set
$$ A_{t}^{F}=\sum_{s \le t}F(X_{s-},X_{s}),\,\, t> 0.$$
We can define a non-local Feynman-Kac semigroup as follows
$$ S_{t}^{F}f(x)=\mathbb{E}_{x}(e^{-A_{t}^{F}}f(X_{t})),$$
where $f$ is a measurable function on $\R^d$. This semigroup was
studied in \cite{S1} and \cite{CS}.

Let $\tilde{F}(x,y)=F(y,x)$, for any $(x,y)\in \R^d \times \R^d$.
In this paper, we always assume both $F$ and $\tilde{F}$ $\in
\bf{J_{d,\alpha}}$.

Recently, sharp two-sided estimates of the density of the
semigroup $\{ S_{t}^{F},t\ge 0\}$ were established in \cite{S2}.
Under the assumption that $X$ is a symmetric $\alpha$-stable-like
process, using a martingale argument and results from \cite{CS},
the following result was established in \cite{S2}: Suppose that
$F\in \bf{J_{d,\alpha}}$ is a symmetric function, the semigroup
$\{ S_{t}^{F},t\ge 0\}$ admits a density $q(t,x,y)$ with respect
to $m$ and that $q$ is jointly continuous on
$(0,\infty)\times\R^d\times\R^d$. Furthermore, there exist
positive constants $C_1,C_2,C_3$ and $C_4$ such that
$$C_{1}e^{-C_{2}t}t^{-\frac{d}{\alpha}}\left(1 \wedge \frac{t^{\frac{1}{\alpha}}}{|x-y|}\right)^{d+\alpha} \leq q(t,x,y) \leq  C_{3}e^{C_{4}t}t^{-\frac{d}{\alpha}}\left(1 \wedge \frac{t^{\frac{1}{\alpha}}}{|x-y|}\right)^{d+\alpha}$$
for all $(t,x,y)\in (0,\infty) \times \R^d \times \R^d$.

The question that we are going to address in this paper is the
following : can we establish the same two-sided estimates for the
density of the Feynman-Kac semigroup of nonsymmetric
$\alpha$-stable-like process $X$ when $F\in \bf{J_{d,\alpha}}$ is
nonsymmetric. The proof of the above result in \cite{S2} can not
be adapted to the case where neither $F$ nor $X$ is symmetric. It
seems that, to answer the question, one has to use some new ideas.
In this paper, we are going to tackle the question above by
combining  the generalization of an idea of \cite{BM}, which was
used to deal with the estimates of the density of continuous
functionals of Brownian motion, with some results on discontinuous
additive functionals.

The content of this paper is organized as follows. In section 2,
we present some preliminary results on discontinuous additive
functionals. In section 3, we establish the two-sided estimates on
the density of Feynman-Kac semigroup under certain assumptions of
$F(x,y)$.

\section{Preliminary Results on Discontinuous Additive Functionals}
For convenience, we use $A_t$ to denote $\sum_{s \le
t}F(X_{s-},X_{s})$ instead of $A^F_t$. We have the following
 formulae for $A^2_{t}$.
$$
    A^2_{t}=2\int_{0}^{t}
    A_{s}\,dA_{s}-\int_{0}^{t}F(X_{s-},X_{s})\,dA_{s},
$$
and
$$
   A^2_{t}=2\int_{0}^{t}(A_{t}- A_{s})\,dA_{s}+\int_{0}^{t}F(X_{s-},X_{s})\,dA_{s}.
$$
The proof is straight forward.

In general, the formulae for $A^n_t$ are given by the following
theorem.
\begin{thm}
\begin{eqnarray*}
            A^n_{t}&=&C^1_{n}\int_{0}^{t}A_{s}^{n-1}\,dA_{s}-C^2_{n}\int_{0}^{t}A_{s}^{n-2}F(X_{s-},X_{s})\,dA_{s}+C^3_{n}\int_{0}^{t}A_{s}^{n-3}F^2(X_{s-},X_{s})\,dA_{s}+\cdots\\
                   &&+(-1)^{i-1}C^{i}_n \int_{0}^{t}A_{s}^{n-i}F^{i-1}(X_{s-},X_{s})\,dA_{s}+\cdots+(-1)^{n-1}C^n_n\int_{0}^{t}F^{n-1}(X_{s-},X_{s})\,dA_{s},\\
            A^n_{t}&=&C^1_{n}\int_{0}^{t}(A_{t}-A_{s})^{n-1}\,dA_{s}+C^2_{n}\int_{0}^{t}(A_{t}-A_{s})^{n-2}F(X_{s-},X_{s})\,dA_{s}\\
                 &&+C^3_{n}\int_{0}^{t}(A_{t}-A_{s})^{n-3}F^2(X_{s-},X_{s})\,dA_{s}+\cdots+C^{i}_n \int_{0}^{t}(A_{t}-A_{s})^{n-i}F^{i-1}(X_{s-},X_{s})\,dA_{s}\\
                 &&+\cdots+C^n_n\int_{0}^{t}F^{n-1}(X_{s-},X_{s})\,dA_{s},\\
          \textrm{ where } C^i_n&=&\frac{n!}{i!(n-i)!}.
\end{eqnarray*}
.
\end{thm}
\pf We use induction to show these two formulae for $A^n_t$ hold
for all $n>1$. It is clear that they are true for $n=2$. Suppose
they hold for $n \leq m-1$, we show they hold for $n=m$.

It follows from the integration by parts formula,
   $$A^m_t=\int_{0}^{t}A_{s-}\,dA^{m-1}_{s}+\int_{0}^{t}A_{s}^{m-1}\,dA_{s}$$
   where
   \begin{eqnarray*}
                &&\int_{0}^{t}A_{s-}\,dA^{m-1}_{s}\\
                &=&\int_{0}^{t}(A_s-F(X_{s-},X_{s}))\,dA^{m-1}_{s}\\
                &=&\int_{0}^{t}A_s\,dA^{m-1}_{s}-\int_{0}^{t}F(X_{s-},X_{s})\,dA^{m-1}_{s}\\
                &=&\int_{0}^{t}A_s(\sum_{i=1}^{m-1}(-1)^{i-1}C_{m-1}^{i}A^{m-1-i}_{s}F^{i-1}(X_{s-},X_{s}))\,dA_{s}\\
                &&-\int_{0}^{t}F(X_{s-},X_{s})(\sum_{j=1}^{m-1}(-1)^{j-1}C_{m-1}^{j}A^{m-1-j}_{s}F^{j-1}(X_{s-},X_{s})\,dA_{s}\\
                && (\textrm{ by the first formula for } A^n_t \textrm{ when } n=m-1 \textrm{ } ) \\
                &=&\sum_{i=1}^{m-1}(-1)^{i-1}C_{m-1}^{i}\int_{0}^{t}A^{m-i}_{s}F^{i-1}(X_{s-},X_{s})\,dA_{s}\\
                &&-\sum_{j=1}^{m-1}(-1)^{j-1}C_{m-1}^{j}\int_{0}^{t}A^{m-1-j}_{s}F^{j}(X_{s-},X_{s})\,dA_{s}\\
                &=&\sum_{i=1}^{m-1}(-1)^{i-1}C_{m-1}^{i}\int_{0}^{t}A^{m-i}_{s}F^{i-1}(X_{s-},X_{s})\,dA_{s}\\
                &&-\sum_{i=2}^{m}(-1)^{i-2}C_{m-1}^{i-1}\int_{0}^{t}A^{m-i}_{s}F^{i-1}(X_{s-},X_{s})\,dA_{s}\\
                &&\textrm{                      } (\textrm{ let }j=i-1 )\\
                &=&\sum_{i=2}^{m-1}(-1)^{i-1}(C_{m-1}^{i}+C_{m-1}^{i-1})\int_{0}^{t}A^{m-i}_{s}F^{i-1}(X_{s-},X_{s})\,dA_{s}\\
                &&+C_{m-1}^{1}\int_{0}^{t}A_{s}^{m-1}\,dA_{s}-(-1)^{m-2}\int_{0}^{t}F^{m-1}(X_{s-},X_{s})\,dA_{s}\\
                &=&\sum_{i=2}^{m-1}(-1)^{i-1}C_{m}^{i}\int_{0}^{t}A^{m-i}_{s}F^{i-1}(X_{s-},X_{s})\,dA_{s}\\
                &&+C_{m-1}^{1}\int_{0}^{t}A_{s}^{m-1}\,dA_{s}-(-1)^{m}\int_{0}^{t}F^{m-1}(X_{s-},X_{s})\,dA_{s}\\
                &&\textrm{   } (\textrm{ by }C_{m-1}^{i}+C_{m-1}^{i-1}=C_{m}^{i}
                ),
\end{eqnarray*}
    thus
$$
A^m_t=\int_{0}^{t}A_{s-}\,dA^{m-1}_{s}+\int_{0}^{t}A_{s}^{m-1}\,dA_{s}=\sum_{i=1}^{m}(-1)^{i-1}C_{m}^{i}\int_{0}^{t}A^{m-i}_{s}F^{i-1}(X_{s-},X_{s})\,dA_{s},
$$
i.e. the first formula for $A^n_t$ holds for n=m.

Now we go to the second formula for $A^n_t$, for $n=m$.
\begin{eqnarray*}
     &&C^1_{m}\int_{0}^{t}(A_{t}-A_{s})^{m-1}\,dA_{s}\\
     &=&C^1_{m}\int_{0}^{t}\sum_{i=0}^{m-1}C_{m-1}^{i}A^{i}_{t}(-1)^{m-1-i}A^{m-1-i}_{s}\,dA_{s}\\
     &=&\sum_{i=0}^{m-1}(-1)^{m-1-i}C^1_{m}C_{m-1}^{i}A^{i}_{t}\int_{0}^{t}A^{m-1-i}_{s}\,dA_{s}\\
      &=&\sum_{i=0}^{m-1}(-1)^{m-1-i}C_{m}^{i}(m-i)A^{i}_{t}\int_{0}^{t}A^{m-1-i}_{s}\,dA_{s}\\
      &&\textrm{        } (\textrm{ by }C_{m}^{1}C_{m-1}^{i}=C_{m}^{i}(m-i)\textrm{ } )\\
      &=&\sum_{i=0}^{m-1}(-1)^{m-1-i}C_{m}^{i}A^{i}_{t}((m-i)\int_{0}^{t}A^{m-1-i}_{s}\,dA_{s})\\
      &=&\sum_{i=0}^{m-1}(-1)^{m-1-i}C_{m}^{i}A^{i}_{t}(A^{m-i}_{t}+\int_{0}^{t}\sum_{k=2}^{m-i}(-1)^{k}C_{m-i}^{k}A^{m-i-k}_{s}F^{k-1}(X_{s-},X_{s})\,dA_{s})\\
      &&\textrm {        }(\textrm{ by the first formula for }A^n_t \textrm{ when } n=m-i \textrm{ } )\\
      &=&\sum_{i=0}^{m-1}(-1)^{m-1-i}C_{m}^{i}A^{i}_{t}A^{m-i}_{t}\\
      &&+\int_{0}^{t}\sum_{i=0}^{m-1}(-1)^{m-1-i}\sum_{k=2}^{m-i}(-1)^{k}C_{m}^{i}C_{m-i}^{k}A^{i}_{t}A^{m-i-k}_{s}F^{k-1}(X_{s-},X_{s})\,dA_{s},\\
\end{eqnarray*}
where
$$\sum_{i=0}^{m-1}(-1)^{m-1-i}C_{m}^{i}A^{i}_{t}A^{m-i}_{t}=(\sum_{i=0}^{m-1}(-1)^{m-1-i}C_{m}^{i})A^{m}_{t}=(1)A^{m}_{t},$$
and
\begin{eqnarray*}
              &&\int_{0}^{t}\sum_{i=0}^{m-1}(-1)^{m-1-i}\sum_{k=2}^{m-i}(-1)^{k}C_{m}^{i}C_{m-i}^{k}A^{i}_{t}A^{m-i-k}_{s}F^{k-1}(X_{s-},X_{s})\,dA_{s}\\
              &=&\int_{0}^{t}\sum_{k=2}^{m}\sum_{i=0}^{m-k}(-1)^{m-k-i-1}C_{m}^{k}C_{m-k}^{i}A^{i}_{t}A^{m-k-i}_{s}F^{k-1}(X_{s-},X_{s})\,dA_{s}\\
              &&\textrm{      }(\textrm{ by }C_{m}^{i}C_{m-i}^{k}=C_{m}^{k}C_{m-k}^{i} \textrm{ and } (-1)^{m-1-i+k}=(-1)^{m-k-i-1}\textrm{ } )\\
              &=&\sum_{k=2}^{m}C_{m}^{k}(-1)^{-1}\int_{0}^{t}(A_{t}-A_{s})^{m-k}F^{k-1}(X_{s-},X_{s})\,dA_{s}\\
              &=&-\sum_{k=2}^{m}C_{m}^{k}\int_{0}^{t}(A_{t}-A_{s})^{m-k}F^{k-1}(X_{s-},X_{s})\,dA_{s},\\
\end{eqnarray*}
therefore
$$C^1_{m}\int_{0}^{t}(A_{t}-A_{s})^{m-1}\,dA_{s}=A^{m}_{t}-\sum_{k=2}^{m}C_{m}^{k}\int_{0}^{t}(A_{t}-A_{s})^{m-k}F^{k-1}(X_{s-},X_{s})\,dA_{s},$$
i.e. the second formula of $A^n_t$ holds for $n=m$. \qed

\section{Density of Feynman-Kac Semigroups Given by Discontinuous Additive Functionals}
From now on we define $q_{0}(t,x,y)=p(t,x,y)$ where $p(t,x,y)$ is
the transition density of $\alpha$-stable-like process $X$ and
satisfies (1.1). By the second formula for $A^n_t$, we have for
any bounded measurable function $g$
\begin{eqnarray*}
       \mathbb{E}_x[A^n_{t}g(X_t)]&=&\sum_{i=1}^nC^{i}_n\mathbb{E}_x[\int_{0}^{t}(A_{t}-A_{s})^{n-i}F^{i-1}(X_{s-},X_{s})g(X_{t})\,dA_{s}]\\
                                  &=&\sum_{i=1}^nC^{i}_n\mathbb{E}_x[\int_{0}^{t}\mathbb{E}_{X_s}\left(A_{t-s}^{n-i}g(X_{t-s})\right)\,d(\sum_{r \le s}F^{i}(X_{r-},X_{r}))]\\
                                  &=&\sum_{i=1}^nC^{i}_n\mathbb{E}_x[\int_{0}^{t}\int_{\R^d}\frac{2C(X_s,y)F^{i}(X_s,y)}{|X_s-y|^{d+\alpha}}\mathbb{E}_{y}\left(A_{t-s}^{n-i}g(X_{t-s})\right)\,m(dy)ds].\\
 \end{eqnarray*}

 We define $q_n(t,x,z)$ as follows,
 $$q_n(t,x,z)=\sum_{i=1}^nC^{i}_n\int_{0}^{t}\int_{\R^d}p(s,x,w)m(dw)\int_{\R^d}\frac{2C(w,y)F^{i}(w,y)}{|w-y|^{d+\alpha}}q_{n-i}(t-s,y,z)\,m(dy)ds.$$
Then by induction, we can show that for any $n>0$
$$\int_{\R^d}q_n(t,x,z)g(z)\,m(dz)=\mathbb{E}_x[A^n_{t}g(X_t)] $$
and
\begin{eqnarray*}
       \mathbb{E}_x[A^n_{t}g(X_t)]&=&\sum_{i=1}^nC^{i}_n\mathbb{E}_x[\int_{0}^{t}\int_{\R^d}\frac{2C(X_s,y)F^{i}(X_s,y)}{|X_s-y|^{d+\alpha}}\int_{\R^d}q_{n-i}(t-s,y,z)g(z)\,m(dz)\,m(dy)ds]\\
                                  &=&\sum_{i=1}^nC^{i}_n\int_{0}^{t}\int_{\R^d}p(s,x,w)m(dw)\int_{\R^d}\frac{2C(w,y)F^{i}(w,y)}{|w-y|^{d+\alpha}}\int_{\R^d}q_{n-i}(t-s,y,z)g(z)\\
                                  &&\cdot m(dz)\,m(dy)ds.\\
                                  \end{eqnarray*}

We assume that there exist positive constants $\overline{C}$, $L$,
$M_0$ and $\overline{M}$ such that $|2C(x,y)|\leq \overline{C}$,
$|F(x,y)|\leq \frac{L}{2}$ and $0 <M_0\leq M(y)\leq \overline{M}$
where $m(dy)=M(y)dy$. Define
$\overline{F}(w,y)=|F(w,y)|+|F(y,w)|$, which is symmetric and
satisfies $|\overline{F}(w,y)| \le L$ . Define
$\overline{p}(t,x,y)=p(t,x,y)+p(t,y,x)$. Then
$\overline{p}(t,x,y)$ is symmetric and satisfies
$$2M_{1}t^{-\frac{d}{\alpha}}\left(1 \wedge \frac{t^{\frac{1}{\alpha}}}{|x-y|}\right)^{d+\alpha}\leq \overline{p}(t,x,y)\leq 2 M_{2}t^{-\frac{d}{\alpha}}\left(1 \wedge \frac{t^{\frac{1}{\alpha}}}{|x-y|}\right)^{d+\alpha},\,\,\forall (t,x,y)\in (0,\infty)\times\R^d\times\R^d. $$

Denote $
(\int_{\R^d}\frac{|\overline{F}(w,y)|}{|w-y|^{d+\alpha}}\,dy)dw$
by $\mu(dw)$ and let
$C_t=\sup_{x\in\R^d}\int_{0}^{t}\int_{\R^d}\overline{p}(s,x,w)\,\mu(dw)ds$.
Then $C_t \downarrow 0$ as $t \downarrow 0$. It is clear that
there exist two positive constants $D_{1}$ and $D_{2}$ such that
$D_{1}\le \int_{\R^d}\overline{p}(t,x,y)\,m(dy)\le D_{2}$, as
$\overline{p}(t,x,y)$ is comparable to $p(t,x,y)$. Let
$\overline{q}_{0}(t,x,z)=\overline{p}(s,x,z)$ and define
$\overline{q}_{n}(t,x,z)$ by
$$
               \overline{q}_n(t,x,z)=\sum_{i=1}^nC^{i}_n\int_{0}^{t}\int_{\R^d}\overline{p}(s,x,w)m(dw)\int_{\R^d}\frac{\overline{C}\overline{F}^{i}(w,y)}{|w-y|^{d+\alpha}}\overline{q}_{n-i}(t-s,y,z)\,m(dy)ds.
$$

We can see that $|q_{n}(t,x,z)| \leq \overline{q}_{n}(t,x,z)$.

Before we move on to the main results, two lemmas are needed.
\begin{lemma}
For any two positive constants $K<1$ and $L$, there exists a
constant $C_{0}(K,L)$ which depends on $K$ and $L$, such that
\begin{equation}
       K^{n-1}+K^{n-2}\frac{L}{2!}+K^{n-3}\frac{L^2}{3!}+\cdots+\cdots+K^{n-i}\frac{L^{i-1}}{i!}+\cdots+\frac{L^{n-1}}{n!}\leq C_{0}(K,L)K^n\,,\\\textrm{ for all }
       n>0.
\end{equation}

\end{lemma}
\pf Use the fact that  there exists $i_{0}\ge0$, such that
$$
      \frac{L^{l-1}}{l!}\leq \left(\frac{K}{2}\right)^{l},
      \,\,\,\textrm{for } l\geq i_{0}.
$$
\qed

\begin{remark}
     For any given $K$ and $L$, we can choose a small $t_{0}$ so that for a given constant $M_1$, $C_tM_1C_{0}(K,L)\le 1$ for $0\leq t \leq t_{0}$.
Thus
\begin{eqnarray*}
      &&C_tM_1\left( K^{n-1}+K^{n-2}\frac{L}{2!}+K^{n-3}\frac{L^2}{3!}+\cdots+\frac{L^{n-1}}{n!}\right)\le C_tM_1C_{0}(K,L)K^n \le K^n.\\
\end{eqnarray*}
\end{remark}

\begin{lemma}
      $\overline{q}_n(t,x,y)$ is symmetric in $x$ and $y$.
\end{lemma}
\pf We know that
\begin{eqnarray*}
\overline{q}_n(t,z,x)&=&\sum_{i_1=1}^nC^{i_1}_n\int_{0}^{t}\int_{\R^d}\overline{p}(s_1,z,w_1)m(dw_1)\int_{\R^d}\frac{\overline{C}\overline{F}^{i_1}(w_1,y_1)}{|w_1-y_1|^{d+\alpha}}\overline{q}_{n-i}(t-s_1,y_1,x)\,m(dy_1)ds_1\\
          &=&\sum_{i_1=1}^nC^{i_1}_n\sum_{i_2=1}^{n-i_1}C^{i_2}_{n-i_1}
          [\int_{0}^{t}\int_{\R^d}
          \overline{p}(s_1,z,w_1)m(dw_1)\int_{\R^d}\frac{\overline{C}\overline{F}^{i_1}(w_1,y_1)}{|w_1-y_1|^{d+\alpha}}\,m(dy_1)ds_1
          \int_{0}^{t-s_1}\\
          &&\cdot\int_{\R^d}\overline{p}(s_2,y_1,w_2)m(dw_2)\int_{\R^d}\frac{\overline{C}\overline{F}^{i_2}(w_2,y_2)}{|w_2-y_2|^{d+\alpha}}\overline{q}_{n-i_1-i_2}(t-s_1-s_2,y_2,x)\,m(dy_2)ds_2]\\
          &\vdots&\\
          &=&\sum_{i_1+i_2+\cdots+i_k=n}C^{i_1}_nC^{i_2}_{n-i_1}\cdot\cdots \cdot C^{i_k}_{n-i_1-i_2-\cdots-i_{k-1}}
          [\int_{0}^{t}\int_{0}^{t-s_1}\cdots\int_{0}^{t-s_1-\cdots-s_{k-1}}\int_{\R^d}\cdots\int_{\R^d}\\
          &&\overline{p}(s_1,z,w_1)\frac{\overline{C}\overline{F}^{i_1}(w_1,y_1)}{|w_1-y_1|^{d+\alpha}}
          \overline{p}(s_2,y_1,w_2)\frac{\overline{C}\overline{F}^{i_2}(w_2,y_2)}{|w_2-y_2|^{d+\alpha}}
          \cdot\cdots \cdot
          \overline{p}(s_k,y_{k-1},w_k)\\
          &&\cdot\frac{\overline{C}\overline{F}^{i_k}(w_k,y_k)}{|w_k-y_k|^{d+\alpha}}\overline{p}(t-s_1-\cdots-s_k,y_k,x)\,ds_1\cdots
            ds_k m(dw_1)\cdots m(dw_k)\\
          &&\cdot m(dy_1)\cdots m(dy_k)]\\
          &=&\sum_{i_1+i_2+\cdots+i_k=n}C^{i_1,\ldots,i_k}_n
          [\int_{0}^{t}\int_{0}^{t-s_1}\cdots\int_{0}^{t-s_1-\cdots-s_{k-1}}\int_{\R^d}\cdots\int_{\R^d}
          \overline{p}(s_1,z,w_1)
          \frac{\overline{C}\overline{F}^{i_1}(w_1,y_1)}{|w_1-y_1|^{d+\alpha}}\\
           && \cdot \overline{p}(s_2,y_1,w_2)\frac{\overline{C}\overline{F}^{i_2}(w_2,y_2)}{|w_2-y_2|^{d+\alpha}}\cdot\cdots \cdot \overline{p}(s_k,y_{k-1},w_k)
          \frac{\overline{C}\overline{F}^{i_k}(w_k,y_k)}{|w_k-y_k|^{d+\alpha}}\\
          && \cdot \overline{p}(t-s_1-\cdots - s_k,y_k,x)\,ds_1\cdots
            ds_k m(dw_1)\cdots m(dw_k) m(dy_1)\cdots m(dy_k)].\\
\end{eqnarray*}
Put $t-s_1-\cdots-s_k=\tilde{s}_1, s_k=\tilde{s}_2,\ldots,
s_l=\tilde{s}_{k+2-l},\ldots, s_2=\tilde{s}_k$. It is easy to see
the absolute value of the Jacobian of this transformation is $1$.
Let $
y_k=\tilde{w}_1,\ldots,y_{l}=\tilde{w}_{k-l+1},\ldots,y_{1}=\tilde{w}_k,\,
w_k=\tilde{y}_1,\ldots,w_{l}=\tilde{y}_{k-l+1},\ldots,w_{1}=\tilde{y}_k$
and $ j_k=i_1,\ldots,j_{l}=i_{k-l+1},\ldots,j_1=i_k $.\\
 Thus the
above equality becomes
\begin{eqnarray*}
    \overline{q}_n(t,z,x)&=&\sum_{j_1+j_2+\cdots+j_k=n}C^{j_1,\ldots,j_k}_n
          [\int_{0}^{t-\tilde{s}_1-\cdots-\tilde{s}_{k-1}}\int_{0}^{t-\tilde{s}_1-\cdots-\tilde{s}_{k-2}}\cdots\int_{0}^{t}\int_{\R^d}\cdots\int_{\R^d}\\
          &&\cdot \overline{p}(t-\tilde{s}_1-\cdots-\tilde{s}_{k},z,\tilde{y}_k)\frac{\overline{C}\overline{F}^{j_k}(\tilde{y}_k,\tilde{w}_k)}{|\tilde{y}_k-\tilde{w}_k|^{d+\alpha}}
            \overline{p}(\tilde{s}_k,\tilde{w}_k,\tilde{y}_{k-1})\\
          &&\cdot\frac{\overline{C}\overline{F}^{j_{k-1}}(\tilde{y}_{k-1},\tilde{w}_{k-1})}{|\tilde{y}_{k-1}-\tilde{w}_{k-1}|^{d+\alpha}}
            \cdot\cdots\cdot
            \overline{p}(\tilde{s}_2,\tilde{w}_2,\tilde{y}_1)
          \frac{\overline{C}\overline{F}^{j_1}(\tilde{y}_1,\tilde{w}_1)}{|\tilde{y}_1-\tilde{w}_1|^{d+\alpha}}\\
          &&\cdot \overline{p}(\tilde{s}_1,\tilde{w}_1,x)\,d\tilde{s}_k\cdots
            d\tilde{s}_1m(d\tilde{y}_k)\cdots m(d\tilde{y}_1)m(d\tilde{w}_k)\cdots
            m(d\tilde{w}_1)].   \end{eqnarray*}
Rearranging the components of the integrand and using the fact
that $\overline{F}(x,y)$ and $\overline{p}(t,x,y)$ are symmetric
in $x$ and $y$, it is easy to see that the above expression for
$\overline{q}_n(t,z,x)$ is equal to $\overline{q}_n(t,x,z)$. \qed

In the proof of the following theorem, we use an idea similar to
that used in \cite{BM} for Brownian motions.
\begin{thm}
There exist two positive constants $ K<1$  and $ M$, and there
exists $ t_{0}> 0$ such that $0<t \leq t_{0}$,
\begin{equation}
      \overline{q}_n(t,x,z)\leq n!MK^nt^{-\frac{d}{\alpha}}\,, \textrm{ for all } n,
\end{equation}
and
\begin{equation}
      \int_{0}^{t}\int_{\R^d}\overline{q}_n(s,x,z)\,\mu(dz)ds\leq C_tn!K^n\,,\textrm{ for all } n.
\end{equation}
\end{thm}
\pf It is clear that when $n=0$, (3.2) and (3.3) hold.
    We assume they hold for $n\leq m-1$, and consider the case $n=m$. Writing $\overline{q}_m(t,x,y)$ into two terms in the following way:
\begin{eqnarray*}
  \overline{q}_m(t,x,z)&=&\sum_{i=1}^mC^{i}_m\int_{0}^{\frac{t}{2}}\int_{\R^d}\overline{p}(s,x,w)m(dw)\int_{\R^d}\frac{\overline{C}\overline{F}^{i}(w,y)}{|w-y|^{d+\alpha}}\overline{q}_{m-i}(t-s,y,z)\,m(dy)ds\\
            &&+\sum_{i=1}^mC^{i}_m\int_{\frac{t}{2}}^{t}\int_{\R^d}\overline{p}(s,x,w)m(dw)\int_{\R^d}\frac{\overline{C}\overline{F}^{i}(w,y)}{|w-y|^{d+{\alpha}}}\overline{q}_{m-i}(t-s,y,z)\,m(dy)ds.
\end{eqnarray*}

Since (3.2) and (3.3) hold for $n\leq m-1$, we have
\begin{eqnarray*}
             &&\sum_{i=1}^mC^{i}_m\int_{0}^{\frac{t}{2}}\int_{\R^d}\overline{p}(s,x,w)m(dw)\int_{\R^d}\frac{\overline{C}\overline{F}^{i}(w,y)}{|w-y|^{d+\alpha}}\overline{q}_{m-i}(t-s,y,z)\,m(dy)ds\\
             &\leq&\sum_{i=1}^mC^{i}_mM{\overline{M}}^2\overline{C}L^{i-1}(m-i)!K^{m-i}\left(\frac{t}{2}\right)^{-\frac{d}{\alpha}}\int_{0}^{\frac{t}{2}}\int_{\R^d}\overline{p}(s,x,w)\,\mu(dw)ds\\
             &\leq&\sum_{i=1}^mC_tC^{i}_mM{\overline{M}}^2\overline{C}L^{i-1}(m-i)!K^{m-i}\left(\frac{t}{2}\right)^{-\frac{d}{\alpha}}.
\end{eqnarray*}
Similarly,
\begin{eqnarray*}
             &&\sum_{i=1}^mC^{i}_m\int_{\frac{t}{2}}^{t}\int_{\R^d}\overline{p}(s,x,w)m(dw)\int_{\R^d}\frac{\overline{C}\overline{F}^{i}(w,y)}{|w-y|^{d+\alpha}}\overline{q}_{m-i}(t-s,y,z)\,m(dy)ds\\
             &\leq&\sum_{i=1}^mC^{i}_mM\left(\frac{t}{2}\right)^{-\frac{d}{\alpha}}\int_{\frac{t}{2}}^{t}\int_{\R^d}\int_{\R^d}\frac{\overline{C}\overline{F}^{i}(w,y)}{|w-y|^{d+\alpha}}\,m(dw)\overline{q}_{m-i}(t-s,y,z)\,m(dy)ds\\
             &\leq&\sum_{i=1}^mC^{i}_mM\left(\frac{t}{2}\right)^{-\frac{d}{\alpha}}\overline{C}L^{i-1}{\overline{M}}^2\int_{\frac{t}{2}}^{t}\int_{\R^d}\overline{q}_{m-i}(t-s,y,z)\,\mu(dy)ds\\
\end{eqnarray*}
\begin{eqnarray*}
             &\leq&\sum_{i=1}^mC^{i}_mM\left(\frac{t}{2}\right)^{-\frac{d}{\alpha}}\overline{C}L^{i-1}{\overline{M}}^2C_t(m-i)!K^{m-i}\\
             &&\textrm{   } (\textrm{ by symmetry, } \overline{q}_{m-i}(t-s,y,z)=\overline{q}_{m-i}(t-s,z,y) \textrm{ and }(3.3) \textrm{ }) \\
             &=&\sum_{i=1}^mC_tC^{i}_mM{\overline{M}}^2\overline{C}L^{i-1}(m-i)!K^{m-i}\left(\frac{t}{2}\right)^{-\frac{d}{\alpha}}.
\end{eqnarray*}
Therefore,
\begin{eqnarray*}
             \overline{q}_m(t,x,z)&\leq& \sum_{i=1}^mC_tC^{i}_m2^{1+{\frac{d}{\alpha}}}M{\overline{M}}^2\overline{C}L^{i-1}(m-i)!K^{m-i}t^{-\frac{d}{\alpha}}\\
                        &=&m!2^{1+{\frac{d}{\alpha}}}M{\overline{M}}^2\overline{C}t^{-\frac{d}{\alpha}}C_t\left(\sum_{i=1}^mK^{m-i}\frac{L^{i-1}}{i!}\right).
\end{eqnarray*}
Let $M_1=2^{1+\frac{d}{\alpha}}{\overline{M}}^2\overline{C}$. Then
by Remark 3.2, we can choose a small $t_{0}$ such that for $0<
t\leq t_{0}$,
$$
              C_tM_1\left(\sum_{i=1}^mK^{m-i}\frac{L^{i-1}}{i!}\right)\leq K^m\,, \textrm{ for any }
              m\,.
$$
Thus
$$
              \overline{q}_m(t,x,y)\leq
              m!MK^mt^{-\frac{d}{\alpha}},
$$
     i.e. (3.2) holds for $n=m$.

Now we show (3.3) holds for $n=m$.
\begin{eqnarray*}
             &&\int_{0}^{t}\int_{\R^d}\overline{q}_m(s,x,z)\,\mu(dz)ds\\
             &=&\int_{0}^{t}\int_{\R^d}\left(\sum_{i=1}^mC^{i}_m\int_{0}^s\int_{\R^d}\overline{p}(u,x,w)m(dw)\int_{\R^d}\frac{\overline{C}\overline{F}^{i}(w,y)}{|w-y|^{d+\alpha}}\overline{q}_{m-i}(s-u,y,z)m(dy)du\right)\mu(dz)ds.
\end{eqnarray*}
Let $s-u=v$, we get
\begin{eqnarray*}
            &&\int_{0}^{t}\int_{\R^d}\overline{q}_m(s,x,z)\,\mu(dz)ds\\
            &=& \int_{0}^{t}\int_{\R^d}\left(\sum_{i=1}^mC^{i}_m\int_{0}^{t-u}\int_{\R^d}\overline{p}(u,x,w)m(dw)\int_{\R^d}\frac{\overline{C}\overline{F}^{i}(w,y)}{|w-y|^{d+\alpha}}\overline{q}_{m-i}(v,y,z)\,\mu(dz)dv\right)\,m(dy)du\\
            &=& \sum_{i=1}^mC^{i}_m\int_{0}^{t}\int_{\R^d}\overline{p}(u,x,w)\left(\int_{0}^{t-u}\int_{\R^d}\overline{q}_{m-i}(v,y,z)\mu(dz)dv\right)\int_{\R^d}\frac{\overline{C}\overline{F}^{i}(w,y)}{|w-y|^{d+\alpha}}\,m(dy)m(dw)du\\
            &\leq& \sum_{i=1}^mC^{i}_m\int_{0}^{t}\int_{\R^d}\overline{p}(u,x,w)C_t(m-i)!K^{m-i}\int_{\R^d}\frac{\overline{C}\overline{F}^{i}(w,y)}{|w-y|^{d+\alpha}}\,m(dy)m(dw)du\\
\end{eqnarray*}
\begin{eqnarray*}
            &=& \sum_{i=1}^mC^{i}_mC_t(m-i)!K^{m-i}\int_{0}^{t}\int_{\R^d}\overline{p}(u,x,w)\int_{\R^d}\frac{\overline{C}\overline{F}^{i}(w,y)}{|w-y|^{d+\alpha}}\,m(dy)m(dw)du\\
            &\leq& \sum_{i=1}^mC^{i}_mC_t(m-i)!K^{m-i}\overline{C}L^{i-1}{\overline{M}}^2\int_{0}^{t}\int_{\R^d}\overline{p}(u,x,w)\,\mu(dw)du\\
            &\leq& C_tC_t{\overline{M}}^2\overline{C}m!\left(\sum_{i=1}^mK^{m-i}\frac{L^{i-1}}{i!}\right).
\end{eqnarray*}
It is clear that for the previous $t_{0}$, when $0<t\leq t_{0}$,
$$
            C_t{\overline{M}}^2\overline{C}\left(\sum_{i=1}^mK^{m-i}\frac{L^{i-1}}{i!}\right)\leq
            K^m.
$$
Thus
$$
            C_tC_t{\overline{M}}^2\overline{C}m!\left(\sum_{i=1}^mK^{m-i}\frac{L^{i-1}}{i!}\right)\leq
            C_tm!K^m,
$$
i.e.
$$
            \int_{0}^{t}\int_{\R^d}\overline{q}_m(s,x,z)\,\mu(dz)ds\leq
            C_tm!K^m.
$$ Therefore (3.3) holds for $n=m$.

\qed

By the above theorem, we have for $0<t\leq t_{0}$,
$$ \sum_{n=0}^{\infty}\frac{\overline{q}_n(t,x,z)}{n!}\leq \sum_{n=0}^{\infty}MK^nt^{-\frac{d}{\alpha}}=M\frac{1}{1-K}t^{-\frac{d}{\alpha}}. $$

Since $|q_n(t,x,z)|\le \overline{q}_n(t,x,z)$,
$\sum_{n=0}^{\infty}\frac{{q}_n(t,x,z)}{n!}$ is uniformly
convergent on $[\epsilon,t_{0}]\times\R^d\times\R^d$, for any
$\epsilon > 0$. Let
$q(t,x,z)=\sum_{n=0}^{\infty}(-1)^n\frac{q_n(t,x,z)}{n!}$. Then
$q(t,x,z)$ is well defined on $(0,t_{0}]\times\R^d\times\R^d$.
Thus we have the following properties for $q(t,x,z)$,
\begin{prop}
 (i) $\int_{\R^d}q(t,x,z)g(z)\,m(dz)=\mathbb{E}_x[e^{-A_{t}}g(X_{t})],$
for any $g$ bounded measurable and any $t>0$, \noindent (ii)
$\int_{\R^d}q(t,x,y)q(s,y,z)m(dy)=q(t+s,x,z),\, \textrm{for any }
t,s
> 0$.
\end{prop}

In the following, we estimate $q(t,x,z)$ from above and from
below.

It is clear that for any $(t,x,y) \in (0,\infty) \times \R^d
\times \R^d$, there exists a positive constant $D_2$ such that
$\int_{R^d}\overline{p}(t,x,y)m(dy)\leq D_2$. For this $D_2$, the
positive constants $M_1$,$\,L$ and $K<1$ given in (1.1), Remark
3.2 and Theorem 3.4, and $\overline{C}$ which is the upper bound
of $|2C(x,y)|$, there exists a large enough positive integer $k$
such that $\frac{L}{(1-\frac{1}{10}-
2^{-\frac{1}{2}})^{d+\alpha}}\frac{1}{2k}\overline{C}D^2_{2}\le
\frac{1}{8}M_{1}.$

Now instead of considering
$$ \overline{q}_1(t,x,z)=\int_{0}^{t}\int_{\R^d}\overline{p}(s,x,w)m(dw)\int_{\R^d}\frac{\overline{C}\overline{F}(w,y)}{|w-y|^{d+\alpha}}\overline{p}(t-s,y,z)\,m(dy)ds,$$
we consider
\begin{eqnarray*}
 \overline{q}_{1,k}(t,x,z)&=&\frac{\overline{q}_1(t,x,z)}{k}\\
                 &=&\int_{0}^{t}\int_{\R^d}\overline{p}(s,x,w)m(dw)\int_{\R^d}\frac{\overline{C}\frac{\overline{F}(w,y)}{k}}{|w-y|^{d+\alpha}}\overline{p}(t-s,y,z)\,m(dy)ds.
\end{eqnarray*}
We have the following theorem
\begin{thm}
 There exists a small $t_{1}$ such that when $0 <t \le t_{1},$
\begin{equation}
 \overline{q}_{1,k}(t,x,z) \le \frac{1}{2}{p}(t,x,z)
,\,\,\forall (x,z)\in \R^d \times \R^d.
\end{equation}
\end{thm}
\pf  We write $\overline{q}_{1,k}(t,x,z)$ into two terms
\begin{eqnarray*}
   \overline{q}_{1,k}(t,x,z)&=&\int_{0}^{\frac{t}{2}}\int_{\R^d}\overline{p}(s,x,w)m(dw)\int_{\R^d}\frac{\overline{C}\frac{\overline{F}(w,y)}{k}}{|w-y|^{d+\alpha}}\overline{p}(t-s,y,z)\,m(dy)ds\\
                   &&+\int_{\frac{t}{2}}^{t}\int_{\R^d}\overline{p}(s,x,w)m(dw)\int_{\R^d}\frac{\overline{C}\frac{\overline{F}(w,y)}{k}}{|w-y|^{d+\alpha}}\overline{p}(t-s,y,z)\,m(dy)ds.
\end{eqnarray*}
First we look at the first term. There are two cases.

Case 1. When $|x-z| \le t^{\frac{1}{\alpha}}$,
\begin{eqnarray*}
    &&\int_{0}^{\frac{t}{2}}\int_{\R^d}\overline{p}(s,x,w)m(dw)\int_{\R^d}\frac{\overline{C}\frac{\overline{F}(w,y)}{k}}{|w-y|^{d+\alpha}}\overline{p}(t-s,y,z)\,m(dy)ds\\
    &\le&\int_{0}^{\frac{t}{2}}\int_{\R^d}\overline{p}(s,x,w)m(dw)\int_{\R^d}\frac{\overline{C}\frac{\overline{F}(w,y)}{k}}{|w-y|^{d+\alpha}}2M_{2}(t-s)^{-\frac{d}{\alpha}}\,m(dy)ds\\
    &\le& 2M_{2}{\overline{M}}^2\overline{C}\frac{1}{k}\int_{0}^{\frac{t}{2}}\int_{\R^d}\overline{p}(s,x,w)(t-s)^{-\frac{d}{\alpha}}\,\mu(dw)ds\\
    &\le& C_{t}M_{2}{\overline{M}}^2\overline{C}\frac{1}{k}2^{1+{\frac{d}{\alpha}}}t^{-\frac{d}{\alpha}}.\\
\end{eqnarray*}
Case 2. When $|x-z| \ge t^{\frac{1}{\alpha}}$. Let
$B_1=\{y\in\R^d|\,|y-z| \ge \frac{1}{10}|x-z|\},
B_2=\{w\in\R^d|\,|w-x| \ge  2^{-\frac{1}{2}}|x-z|\}$ and
$B_3=\{(w,y)\in \R^d \times \R^d|\,|y-z| <
\frac{1}{10}|x-z|,\,|w-x| <  2^{-\frac{1}{2}}|x-z|\} $. On $B_3$,
we have $|w-y| \ge (1-\frac{1}{10}- 2^{-\frac{1}{2}} )|x-z|$.
\begin{eqnarray*}
            &&\int_{0}^{\frac{t}{2}}\int_{\R^d}\overline{p}(s,x,w)m(dw)\int_{\R^d}\frac{\overline{C}\frac{\overline{F}(w,y)}{k}}{|w-y|^{d+\alpha}}\overline{p}(t-s,y,z)\,m(dy)ds\\
            &\leq&\int_{0}^{\frac{t}{2}}\int_{\R^d}\overline{p}(s,x,w)m(dw)\int_{\R^d}\frac{\overline{C}\frac{\overline{F}(w,y)}{k}}{|w-y|^{d+\alpha}}\overline{p}(t-s,y,z)1_{B_1}(y)\,m(dy)ds\\
                            &&+\int_{0}^{\frac{t}{2}}\int_{\R^d}\overline{p}(s,x,w)m(dw)\int_{\R^d}\frac{\overline{C}\frac{\overline{F}(w,y)}{k}}{|w-y|^{d+\alpha}}\overline{p}(t-s,y,z)\,1_{B_2}(w)\,m(dy)ds\\
                            &&+\int_{0}^{\frac{t}{2}}\int_{\R^d}\overline{p}(s,x,w)m(dw)\int_{\R^d}\frac{\overline{C}\frac{\overline{F}(w,y)}{k}}{|w-y|^{d+\alpha}}\overline{p}(t-s,y,z)\,1_{B_3}(w,y)\,m(dy)ds\\
            &\leq&\int_{0}^{\frac{t}{2}}\int_{\R^d}\overline{p}(s,x,w)m(dw)\int_{\R^d}\frac{\overline{C}\frac{\overline{F}(w,y)}{k}}{|w-y|^{d+\alpha}}2M_2\overline{M}10^{d+\alpha}\frac{(t-s)}{|x-z|^{d+\alpha}}1_{B_1}(y)\,dyds\\
                            &&+\int_{0}^{\frac{t}{2}}\int_{\R^d}2M_2{\overline{M}}^22^{\frac{1}{2}(d+\alpha)}\frac{s}{{|x-z|}^{d+\alpha}}dw\int_{\R^d}\frac{\overline{C}\frac{\overline{F}(w,y)}{k}}{|w-y|^{d+\alpha}}\overline{p}(t-s,y,z)\,1_{B_2}(w)\,dyds\\
                            &&+\frac{L}{(1-\frac{1}{10}- 2^{-\frac{1}{2}})^{d+\alpha}}\frac{1}{k}\overline{C}\int_{0}^{\frac{t}{2}}\int_{\R^d}\overline{p}(s,x,w)m(dw)\int_{\R^d}\frac{1}{{|x-z|}^{d+\alpha}}\overline{p}(t-s,y,z)\,1_{B_3}(w,y)\,\\
                            &&\cdot m(dy)ds\\
                            &\leq&C_t2M_{2}{\overline{M}}^2\overline{C}\frac{1}{k}10^{d+\alpha}\frac{t}{|x-z|^{d+\alpha}}\\
                            &&+C_t2M_{2}{\overline{M}}^2\overline{C}\frac{1}{k}2^{\frac{1}{2}(d+\alpha)}\frac{\frac{t}{2}}{|x-z|^{d+\alpha}}\\
                            &&+\frac{L}{(1-\frac{1}{10}- 2^{-\frac{1}{2}})^{d+\alpha}}\frac{1}{k}\overline{C}D^{2}_{2}\frac{\frac{t}{2}}{|x-z|^{d+\alpha}}\\
                            &&( \textrm{ by } \int_{\R^d}\overline{p}(s,x,w)m(dw)\le D_{2} \textrm{ and }\int_{\R^d}\overline{p}(t-s,y,z)m(dy)\le D_{2} \textrm{ } )\\
            &\leq&C_tM_{2}{\overline{M}}^2\overline{C}\frac{1}{k}(2\cdot 10^{d+\alpha}+2^{\frac{1}{2}(d+\alpha)})\frac{t}{|x-z|^{d+\alpha}}+\frac{1}{8}p(t,x,z).\\
\end{eqnarray*}
Since $C_t \downarrow 0$ as $t \downarrow 0$, then for both case 1 and case 2, we can find a small $t_{11}$ such that when $0 < t \leq t_{11}$,
$$\int_{0}^{\frac{t}{2}}\int_{\R^d}\overline{p}(s,x,w)m(dw)\int_{\R^d}\frac{\overline{C}\frac{\overline{F}(w,y)}{k}}{|w-y|^{d+\alpha}}\overline{p}(t-s,y,z)\,m(dy)ds \leq \frac{1}{4}p(t,x,z).$$

For the second term of $\overline{q}_{1,k}(t,x,z)$:
$$\int_{\frac{t}{2}}^{t}\int_{\R^d}\overline{p}(s,x,w)m(dw)\int_{\R^d}\frac{\overline{C}\frac{\overline{F}(w,y)}{k}}{|w-y|^{d+\alpha}}\overline{p}(t-s,y,z)\,m(dy)ds.$$
Letting $t-s=\tilde{s}$, the second term becomes
$$\int_{0}^{\frac{t}{2}}\int_{\R^d}\overline{p}(t-\tilde{s},x,w)m(dw)\int_{\R^d}\frac{\overline{C}\frac{\overline{F}(w,y)}{k}}{|w-y|^{d+\alpha}}\overline{p}(\tilde{s},y,z)\,m(dy)d\tilde{s}.$$
There are two cases.

Case a. When $|x-z| \le t^{\frac{1}{\alpha}}$,
\begin{eqnarray*}
    &&\int_{0}^{\frac{t}{2}}\int_{\R^d}\overline{p}(t-\tilde{s},x,w)m(dw)\int_{\R^d}\frac{\overline{C}\frac{\overline{F}(w,y)}{k}}{|w-y|^{d+\alpha}}\overline{p}(\tilde{s},y,z)\,m(dy)d\tilde{s}\\
    &\le&\int_{0}^{\frac{t}{2}}\int_{\R^d}2M_{2}(t-\tilde{s})^{-\frac{d}{\alpha}}m(dw)\int_{\R^d}\frac{\overline{C}\frac{\overline{F}(w,y)}{k}}{|w-y|^{d+\alpha}}\overline{p}(\tilde{s},y,z)\,m(dy)d\tilde{s}\\
    &\le&M_{2}{2\overline{M}}^2\overline{C}\frac{1}{k}\int_{0}^{\frac{t}{2}}\int_{\R^d}\overline{p}(\tilde{s},y,z)(t-\tilde{s})^{-\frac{d}{\alpha}}\,\mu(dy)ds\\
    &\le& C_{t}2M_{2}{\overline{M}}^2\overline{C}\frac{1}{k}2^{\frac{d}{\alpha}}t^{-\frac{d}{\alpha}}\\
    && (\textrm{ by symmetry of } \overline{p}(\tilde{s},y,z) ).
\end{eqnarray*}

Case b. When $|x-z| \ge t^{\frac{1}{\alpha}}$. Let
$\tilde{B_1}=\{y\in\R^d|\,|y-z| \ge \frac{1}{10}|x-z|\},
\tilde{B_2}=\{w\in\R^d|\,|w-x| \ge  2^{-\frac{1}{2}}|x-z|\}$ and
$\tilde{B_3}=\{(w,y)\in \R^d \times \R^d|\,|y-z| <
\frac{1}{10}|x-z|,\,|w-x| <  2^{-\frac{1}{2}}|x-z|\} $. On
$\tilde{B_3}$, we have $|w-y| \ge (1-\frac{1}{10}-
2^{-\frac{1}{2}} )|x-z|$.
\begin{eqnarray*}
            &&\int_{0}^{\frac{t}{2}}\int_{\R^d}\overline{p}(t-\tilde{s},x,w)m(dw)\int_{\R^d}\frac{\overline{C}\frac{\overline{F}(w,y)}{k}}{|w-y|^{d+\alpha}}\overline{p}(\tilde{s},y,z)\,m(dy)d\tilde{s}\\
            &\leq&\int_{0}^{\frac{t}{2}}\int_{\R^d}\overline{p}(t-\tilde{s},x,w)m(dw)\int_{\R^d}\frac{\overline{C}\frac{\overline{F}(w,y)}{k}}{|w-y|^{d+\alpha}}\overline{p}(\tilde{s},y,z)1_{\tilde{B_1}}(y)\,m(dy)d\tilde{s}\\
                            &&+\int_{0}^{\frac{t}{2}}\int_{\R^d}\overline{p}(t-\tilde{s},x,w)m(dw)\int_{\R^d}\frac{\overline{C}\frac{\overline{F}(w,y)}{k}}{|w-y|^{d+\alpha}}\overline{p}(\tilde{s},y,z)\,1_{\tilde{B_2}}(w)\,m(dy)d\tilde{s}\\
                            &&+\int_{0}^{\frac{t}{2}}\int_{\R^d}\overline{p}(t-\tilde{s},x,w)m(dw)\int_{\R^d}\frac{\overline{C}\frac{\overline{F}(w,y)}{k}}{|w-y|^{d+\alpha}}\overline{p}(\tilde{s},y,z)\,1_{\tilde{B_3}}(w,y)\,m(dy)d\tilde{s}\\
            &\leq&\int_{0}^{\frac{t}{2}}\int_{\R^d}\overline{p}(t-\tilde{s},x,w)m(dw)\int_{\R^d}\frac{\overline{C}\frac{\overline{F}(w,y)}{k}}{|w-y|^{d+\alpha}}2M_2\overline{M}10^{d+\alpha}\frac{\tilde{s}}{|x-z|^{d+\alpha}}1_{\tilde{B_1}}(y)\,dyd\tilde{s}\\
                            &&+\int_{0}^{\frac{t}{2}}\int_{\R^d}2M_2{\overline{M}}^22^{\frac{1}{2}(d+\alpha)}\frac{(t-\tilde{s})}{{|x-z|}^{d+\alpha}}dw\int_{\R^d}\frac{\overline{C}\frac{\overline{F}(w,y)}{k}}{|w-y|^{d+\alpha}}\overline{p}(\tilde{s},y,z)\,1_{B_2}(w)\,dyd\tilde{s}\\
                            &&+\frac{L}{(1-\frac{1}{10}- 2^{-\frac{1}{2}})^{d+\alpha}}\frac{1}{k}\overline{C}\int_{0}^{\frac{t}{2}}\int_{\R^d}\overline{p}(t-\tilde{s},x,w)m(dw)\int_{\R^d}\frac{1}{{|x-z|}^{d+\alpha}}\overline{p}(\tilde{s},y,z)\,1_{B_3}(w,y)\,\\
                            &&\cdot m(dy)d\tilde{s}\\
\end{eqnarray*}
\begin{eqnarray*}
            &\leq&C_t2M_{2}{\overline{M}}^2\overline{C}\frac{1}{k}10^{d+\alpha}\frac{t}{|x-z|^{d+\alpha}}\\
                            &&+C_t2M_{2}{\overline{M}}^2\overline{C}\frac{1}{k}2^{\frac{1}{2}(d+\alpha)}\frac{t}{|x-z|^{d+\alpha}}\\
                            &&+\frac{L}{(1-\frac{1}{10}- 2^{-\frac{1}{2}})^{d+\alpha}}\frac{1}{k}\overline{C}D^2_{2}\frac{\frac{t}{2}}{|x-z|^{d+\alpha}}\\
                            &&( \textrm{ by } \int_{\R^d}\overline{p}(t-\tilde{s},x,w)m(dw)\le D_{2} \textrm{ and }\int_{\R^d}\overline{p}(\tilde{s},y,z)m(dy)\le D_{2} \textrm{ } )\\
            &\leq&C_t2M_{2}{\overline{M}}^2\overline{C}\frac{1}{k}(10^{d+\alpha}+2^{\frac{1}{2}(d+\alpha)})\frac{t}{|x-z|^{d+\alpha}}+\frac{1}{8}p(t,x,z).
\end{eqnarray*}
Since $C_t \downarrow 0$ as $t \downarrow 0$, then for both case a
and case b, we can find a small $t_{12}$ such that when $0 < t
\leq t_{12}$,
 $$\int_{0}^{\frac{t}{2}}\int_{\R^d}\overline{p}(t-\tilde{s},x,w)m(dw)\int_{\R^d}\frac{\overline{C}\frac{\overline{F}(w,y)}{k}}{|w-y|^{d+\alpha}}\overline{p}(\tilde{s},y,z)\,m(dy)d\tilde{s}\leq \frac{1}{4}p(t,x,z),$$
i.e.
 $$\int_{\frac{t}{2}}^{t}\int_{\R^d}\overline{p}(s,x,w)m(dw)\int_{\R^d}\frac{\overline{C}\frac{\overline{F}(w,y)}{k}}{|w-y|^{d+\alpha}}\overline{p}(t-s,y,z)\,m(dy)ds \leq \frac{1}{4}p(t,x,z).$$
Let $t_{1}=\textrm{min}(t_{11},t_{12})$. Then when  $0 < t \leq
t_{1},$
$$ \overline{q}_{1,k}(t,x,z) \le \frac{1}{2}p(t,x,z),\,\, \forall (x,z)\in \R^d \times \R^d. $$
\qed

Define $q_{1,k}(t,x,z)=\frac{q_{1}(t,x,z)}{k}$. It  is clear that
$|q_{1,k}(t,x,z)|\leq \overline{q}_{1,k}(t,x,z)$.  Theorem 3.6
implies that when $0 < t \leq t_{1},$
$$ p(t,x,z)- q_{1,k}(t,x,z) \ge p(t,x,z)- \overline{q}_{1,k}(t,x,z) \ge \frac{1}{2}p(t,x,z).$$
We know
$\int_{\R^d}q_{1,k}(t,x,z)g(z)\,m(dz)=\mathbb{E}_x[\frac{A_{t}}{k}g(X_t)]$,
for any $g$ measurable. Since $1-\frac{A_{t}}{k} \leq
e^{-\frac{A_{t}}{k}}$, we have
$$
 \frac{1}{|B_{r}|}\mathbb{E}_x[(1-\frac{A_{t}}{k})1_{B_{r}}(X_{t})] \leq \frac{1}{|B_{r}|}\mathbb{E}_x[e^{-\frac{A_{t}}{k}}1_{B_{r}}(X_{t})].\\
$$
Thus
\begin{eqnarray*}
    \frac{1}{2}\frac{1}{|B_{r}|}\mathbb{E}_x[1_{B_{r}}(X_{t})] &\leq& \frac{1}{|B_{r}|}\mathbb{E}_x[e^{-\frac{A_{t}}{k}}1_{B_{r}}(X_{t})]\\
                                                      &\leq&(\frac{1}{|B_{r}|}\mathbb{E}_x[e^{-A_{t}}1_{B_{r}}(X_{t})])^{\frac{1}{k}}(\frac{1}{|B_{r}|}\mathbb{E}_x[1_{B_{r}}(X_{t})])^{1-\frac{1}{k}}\\
&&(\textrm{ by H\"older inequality }).\\
\end{eqnarray*}
Therefore
\begin{eqnarray*}
\frac{
\frac{1}{2}\frac{1}{|B_{r}|}\mathbb{E}_x[1_{B_{r}}(X_{t})]}{(\frac{1}{|B_{r}|}\mathbb{E}_x[1_{B_{r}}(X_{t})])^{1-\frac{1}{k}}}
&\leq& (\frac{1}{|B_{r}|}\mathbb{E}_x[e^{-A_{t}}1_{B_{r}}(X_{t})])^{\frac{1}{k}},\\
\textrm{ i.e. }\\
    \frac{1}{2}(\frac{1}{|B_{r}|}\mathbb{E}_x[1_{B_{r}}(X_{t})])^{\frac{1}{k}} &\leq& (\frac{1}{|B_{r}|}\mathbb{E}_x[e^{-A_{t}}1_{B_{r}}(X_{t})])^{\frac{1}{k}},\\
\textrm{ i.e. }\\
    \frac{1}{2^k}\frac{1}{|B_{r}|}\mathbb{E}_x[1_{B_{r}}(X_{t})] &\leq&
    \frac{1}{|B_{r}|}\mathbb{E}_x[e^{-A_{t}}1_{B_{r}}(X_{t})].
\end{eqnarray*}
Let $r\downarrow 0$, we have
$$ \frac{1}{2^k}p(t,x,z) \leq q(t,x,z). $$
Therefore when $0 < t \leq t_{0},$
$$ \frac{1}{2^k}M_{1}t^{-\frac{d}{\alpha}}\left(1 \wedge \frac{t^{\frac{1}{\alpha}}}{|x-z|}\right)^{d+\alpha} \leq q(t,x,z). $$
Applying (iv) of Proposition 3.5, we have
$$q(t,x,y)\ge {C}_{3}e^{-{C}_{4}t}t^{-\frac{d}{\alpha}}\left(1 \wedge \frac{t^{\frac{1}{\alpha}}}{|x-y|}\right)^{d+\alpha},\,\,\forall (t,x,y)\in (0,\infty) \times \R^d \times \R^d,$$
where $C_{3}$ and $C_{4}$ are positive constants.

Next we establish the upper bound.

It is clear that for  the positive constants $L$ and $K<1$ given
in Remark 3.2 and Theorem 3.4, and
$\overline{M}$,$\,\overline{C}$, which are the upper bounds for
$M(y)$ and $|2C(x,y)|$ respectively, there exists a constant $
\tilde{C} \ge 1$ such that
$$ L^{n-1}{\overline{M}}^2\overline{C} \le
 \tilde{C}\frac{1}{2}n!K^{n}, \,\, \forall n \ge 1.$$
Suppose that $g \ge 0$ is a measurable function and $g \le
C_{g}\min(\frac{1}{D_{2}},1)$, where $C_{g} \ge 1$ is a constant,
then we have the following
\begin{prop}
    There exists $t_{2}\ge 0$ such that when $0< t \le t_{2}$,
\begin{equation}
    \int_{\R^d}\overline{q}_n(t,x,z)g(z)\,m(dz)
    \le\tilde{C}C_{g}C_{t}n!K^{n},\,\, \forall n \ge 1.
\end{equation}
\end{prop}
\pf When $n=1$,
\begin{eqnarray*}
    &&\int_{\R^d}\overline{q}_1(t,x,z)g(z)\,m(dz)\\
      &=&\int_{\R^d}\int_{0}^{t}\int_{\R^d}\overline{p}(s,x,w)m(dw)\int_{\R^d}\frac{\overline{C}\overline{F}(w,y)}{|w-y|^{d+\alpha}}\overline{p}(t-s,y,z)\,m(dy)ds g(z)\,m(dz)\\
\end{eqnarray*}
\begin{eqnarray*}
      &\le&\int_{0}^{t}\int_{\R^d}\overline{p}(s,x,w)m(dw)\int_{\R^d}\frac{\overline{C}\overline{F}(w,y)}{|w-y|^{d+\alpha}}C_{g}\overline{M}\,dy\,ds\\
      &&(\textrm{ by }\int_{\R^d}\overline{p}(t-s,y,z)g(z)\,m(dz) \le C_{g})\\
      &\le& {\overline{M}}^2\overline{C}C_{g}C_{t}\le \tilde{C}\frac{1}{2}KC_{g}C_{t}\le \tilde{C}C_{g}C_{t}K.\\
\end{eqnarray*}
Thus (3.5) holds for $n=1$.

Suppose it holds for $n \le m-1$, we show that it holds for
    $n=m$.
\begin{eqnarray*}
      &&\int_{\R^d}\overline{q}_m(t,x,z)g(z)\,m(dz)\\
      &=&\int_{\R^d}\sum_{i=1}^mC^{i}_m\int_{0}^{t}\int_{\R^d}\overline{p}(s,x,w)m(dw)\int_{\R^d}\frac{\overline{C}\overline{F}^{i}(w,y)}{|w-y|^{d+\alpha}}\overline{q}_{m-i}(t-s,y,z)\,m(dy)ds g(z)\,m(dz)\\
      &=&\sum_{i=1}^mC^{i}_m\int_{0}^{t}\int_{\R^d}\overline{p}(s,x,w)m(dw)\int_{\R^d}\frac{\overline{C}\overline{F}^{i}(w,y)}{|w-y|^{d+\alpha}}\int_{\R^d}\overline{q}_{m-i}(t-s,y,z)g(z)\,m(dz)\,m(dy)ds\\
      &\leq&\sum_{i=1}^{m-1}C^{i}_m\int_{0}^{t}\int_{\R^d}\overline{p}(s,x,w)m(dw)\int_{\R^d}\frac{\overline{C}\overline{F}^{i}(w,y)}{|w-y|^{d+\alpha}}\,m(dy)ds\tilde{C}C_{g}C_{t}(m-i)!K^{m-i}\\
      &&+\int_{0}^{t}\int_{\R^d}\overline{p}(s,x,w)m(dw)\int_{\R^d}\frac{\overline{C}\overline{F}^{m}(w,y)}{|w-y|^{d+\alpha}}\,m(dy)dsC_{g}\\
      &=&\sum_{i=1}^{m-1}C^{i}_mC_{t}L^{i-1} {\overline{M}}^2\overline{C}\tilde{C}C_{g}C_{t}(m-i)!K^{m-i}+L^{m-1}{\overline{M}}^2\overline{C}C_{g}C_{t}\\
      &=&m!\sum_{i=1}^{m-1}(\frac{L^{i-1}K^{m-i}}{i!})C_{t}{\overline{M}}^2\overline{C}\tilde{C}C_{g}C_{t}+L^{m-1}{\overline{M}}^2\overline{C}C_{g}C_{t}.\\
\end{eqnarray*}
Since $C_{t} \downarrow 0 $ as $ t \downarrow 0$, $\exists t_{2}
\ge 0 $ such that when $0< t \le t_{2}$
$$\sum_{i=1}^{n-1}(\frac{L^{i-1}K^{n-i}}{i!})C_{t}{\overline{M}}^2\overline{C}\le \frac{1}{2}K^{n},\,\,\forall n \ge 2,$$
 by the choice of $\tilde{C}$,
$$L^{n-1}{\overline{M}}^2\overline{C} \le
 \tilde{C}\frac{1}{2}n!K^{n},\,\,\forall n \ge 1,$$
 Thus
\begin{eqnarray*}
 \int_{\R^d}\overline{q}_m(t,x,z)g(z)\,dz
 &\le&\frac{1}{2}m!K^{m}\tilde{C}C_{g}C_{t}+\frac{1}{2}m!K^{m}\tilde{C}C_{g}C_{t}\\
                             &=&\tilde{C}C_{g}C_{t}m!K^{m},
\end{eqnarray*}
i.e. the statement holds for $n=m$. \qed

For the $L$, $K$, $\overline{C}$ and $D_2$ given above, it is
clear that there exists $\tilde{C_{2}} \ge 1$ such that
$$\frac{L^n}{(1-\frac{1}{10}-2^{-\frac{1}{2}})^{d+\alpha}}\frac{\overline{C}D^2_{2}}{2} \le
       \frac{1}{8}\tilde{C_{2}}n!K^{n},\,\,\forall n \ge 0. $$
We claim that
\begin{thm}
      There exist $ t_{3} >0 \textrm{ and }\tilde{C_{1}} \ge 1$ such
      that when $0 < t \le t_{3}$,
\begin{equation}
 \overline{q}_{n}(t,x,z) \le
\tilde{C_{1}}n!K^{n}t^{-\frac{d}{\alpha}}\left(1 \wedge
\frac{t^{\frac{1}{\alpha}}}{|x-z|}\right)^{d+\alpha},\,\,\forall n
\ge 0.
\end{equation}

\end{thm}
\pf Since $\overline{q}_{0}(t,x,z)=\overline{p}(t,x,z)$, there
exist $ t_{13}>0 \textrm{ and } \tilde{C_{1}} \ge \tilde{C_{2}}$
such that when $0 < t \le t_{13}$
       $$\overline{q}_{0}(t,x,z) \le \tilde{C_{1}}t^{-\frac{d}{\alpha}}\left(1 \wedge \frac{t^{\frac{1}{\alpha}}}{|x-z|}\right)^{d+\alpha},$$
i.e. the statement holds for $n=0$.
    Suppose it is true for $n \le m-1$. We show that it holds for
    $n=m$.
We write $\overline{q}_{m}(t,x,z)$ into two terms
\begin{eqnarray*}
\overline{q}_m(t,x,z)&=&\sum_{i=1}^mC^{i}_m\int_{0}^{\frac{t}{2}}\int_{\R^d}\overline{p}(s,x,w)m(dw)\int_{\R^d}\frac{\overline{C}\overline{F}^{i}(w,y)}{|w-y|^{d+\alpha}}\overline{q}_{m-i}(t-s,y,z)\,m(dy)ds\\
           &&+\sum_{i=1}^mC^{i}_m\int_{\frac{t}{2}}^{t}\int_{\R^d}\overline{p}(s,x,w)m(dw)\int_{\R^d}\frac{\overline{C}\overline{F}^{i}(w,y)}{|w-y|^{d+\alpha}}\overline{q}_{m-i}(t-s,y,z)\,m(dy)ds.\\
\end{eqnarray*}
First we look at the first term. There are two cases:

 Case 1. When $|x-z| \le t^{\frac{1}{\alpha}}$,
\begin{eqnarray*}
        &&\sum_{i=1}^mC^{i}_m\int_{0}^{\frac{t}{2}}\int_{\R^d}\overline{p}(s,x,w)m(dw)\int_{\R^d}\frac{\overline{C}\overline{F}^{i}(w,y)}{|w-y|^{d+\alpha}}\overline{q}_{m-i}(t-s,y,z)\,m(dy)ds\\
        &\le& \sum_{i=1}^mC^{i}_m\int_{0}^{\frac{t}{2}}\int_{\R^d}\overline{p}(s,x,w)m(dw)\int_{\R^d}\frac{\overline{C}\overline{F}(w,y)}{|w-y|^{d+\alpha}}\,dydsL^{i-1}\overline{M}\tilde{C_{1}}(m-i)!K^{m-i}(\frac{t}{2})^{-\frac{d}{\alpha}}\\
        &\le&
        \sum_{i=1}^mC^{i}_mC_{t}{\overline{M}}^2\overline{C}L^{i-1}\tilde{C_{1}}(m-i)!K^{m-i}(\frac{t}{2})^{-\frac{d}{\alpha}}\\
        &=&m!\sum_{i=1}^m(\frac{L^{i-1}K^{m-i}}{i!}){\overline{M}}^2\overline{C}(\frac{1}{2})^{-\frac{d}{\alpha}}C_{t}\tilde{C_{1}}t^{-\frac{d}{\alpha}}.
\end{eqnarray*}
  Since there exists $ t_{23} >0 \textrm{ and } t_{23} \le t_{13} $ such that when $0 < t \le
  t_{23}$,
       $$\sum_{i=1}^n(\frac{L^{i-1}K^{n-i}}{i!}){\overline{M}}^2\overline{C}(\frac{1}{2})^{-\frac{d}{\alpha}}C_{t}\le
       \frac{1}{2}K^{n},\,\,\forall n \ge 1,$$
  thus in case 1, when $0 < t \le t_{23}$,
\begin{eqnarray*}
\sum_{i=1}^mC^{i}_m\int_{0}^{\frac{t}{2}}\int_{\R^d}\overline{p}(s,x,w)m(dw)\int_{\R^d}\frac{\overline{C}\overline{F}^{i}(w,y)}{|w-y|^{d+\alpha}}\overline{q}_{m-i}(t-s,y,z)\,m(dy)ds
        \le \frac{1}{2}\tilde{C_{1}}m!K^{m}t^{-\frac{d}{\alpha}}.
\end{eqnarray*}

Case 2. When $|x-z|\ge t^{\frac{1}{\alpha}}$. Let
$B_1=\{y\in\R^d|\,|y-z| \ge \frac{1}{10}|x-z|\},
B_2=\{w\in\R^d|\,|w-x| \ge  2^{-\frac{1}{2}}|x-z|\}$ and
$B_3=\{(w,y)\in \R^d \times \R^d|\,|y-z| <  \frac{1}{10}|x-z|,\,
|w-x| <  2^{-\frac{1}{2}}|x-z|\} $. On $B_3$, we have $|w-y| \ge
(1-\frac{1}{10}- 2^{-\frac{1}{2}}
  )|x-z|$.
\begin{eqnarray*}
       &&\sum_{i=1}^mC^{i}_m\int_{0}^{\frac{t}{2}}\int_{\R^d}\overline{p}(s,x,w)m(dw)\int_{\R^d}\frac{\overline{C}\overline{F}^{i}(w,y)}{|w-y|^{d+\alpha}}\overline{q}_{m-i}(t-s,y,z)\,m(dy)ds\\
       &\le&\sum_{i=1}^mC^{i}_m\int_{0}^{\frac{t}{2}}\int_{\R^d}\overline{p}(s,x,w)m(dw)\int_{\R^d}\frac{\overline{C}\overline{F}^{i}(w,y)}{|w-y|^{d+\alpha}}\overline{q}_{m-i}(t-s,y,z)1_{B_{1}}(y)\,m(dy)ds\\
       &&+\sum_{i=1}^mC^{i}_m\int_{0}^{\frac{t}{2}}\int_{\R^d}\overline{p}(s,x,w)m(dw)\int_{\R^d}\frac{\overline{C}\overline{F}^{i}(w,y)}{|w-y|^{d+\alpha}}\overline{q}_{m-i}(t-s,y,z)1_{B_{2}}(w)  \,m(dy)ds\\
       &&+\sum_{i=1}^mC^{i}_m\int_{0}^{\frac{t}{2}}\int_{\R^d}\overline{p}(s,x,w)m(dw)\int_{\R^d}\frac{\overline{C}\overline{F}^{i}(w,y)}{|w-y|^{d+\alpha}}\overline{q}_{m-i}(t-s,y,z)1_{B_{3}}(w,y) \,m(dy)ds\\
       &\le& \sum_{i=1}^mC^{i}_m\int_{0}^{\frac{t}{2}}\int_{\R^d}\overline{p}(s,x,w)m(dw)\int_{\R^d}\frac{\overline{C}\overline{F}(w,y)}{|w-y|^{d+\alpha}} \overline{M}L^{i-1}\tilde{C_{1}}(m-i)!K^{m-i}10^{d+\alpha}\frac{(t-s)}{|x-z|^{d+\alpha}}\,dyds\\
       &&+\sum_{i=1}^mC^{i}_m\int_{0}^{\frac{t}{2}}\int_{\R^d}\,m(dw)\tilde{C_{1}}{\overline{M}}^22^{\frac{1}{2}(d+\alpha)}\frac{s}{{|x-z|}^{d+\alpha}}\int_{\R^d}\frac{\overline{C}\overline{F}(w,y)}{|w-y|^{d+\alpha}}\overline{q}_{m-i}(t-s,y,z)\,dydsL^{i-1}\\
       &&+\sum_{i=1}^{m-1}C^{i}_m\frac{L^i}{(1-\frac{1}{10}- 2^{-\frac{1}{2}})^{d+\alpha}}\overline{C}\int_{0}^{\frac{t}{2}}\int_{\R^d}\overline{p}(s,x,w)m(dw)\int_{\R^d}\frac{1}{|x-z|^{d+\alpha}}\overline{q}_{m-i}(t-s,y,z)\,\\
       &&\cdot m(dy)ds+\frac{L^m}{(1-\frac{1}{10}-2^{-\frac{1}{2}})^{d+\alpha}}\overline{C}D^2_{2}\frac{\frac{t}{2}}{|x-z|^{d+\alpha}}\\
       &\le& \sum_{i=1}^mC^{i}_mC_{t}
       {\overline{M}}^2\overline{C}L^{i-1}\tilde{C_{1}}(m-i)!K^{m-i}10^{d+\alpha}\frac{t}{|x-z|^{d+\alpha}}\\
       &&+\sum_{i=1}^mC^{i}_m\tilde{C_{1}}{\overline{M}}^2\overline{C}2^{\frac{1}{2}(d+\alpha)}\frac{t}{{|x-z|}^{d+\alpha}}C_{t}(m-i)!K^{m-i}L^{i-1}\\
       &&( \textrm{ by symmetry of } \overline{q}_{m-i}(t-s,y,z)\textrm{ and Theorem 3.4 } )\\
       &&+\sum_{i=1}^{m-1}C^{i}_m\frac{L^i}{(1-\frac{1}{10}- 2^{-\frac{1}{2}})^{d+\alpha}}\overline{C}D_{2}\frac{t}{|x-z|^{d+\alpha}}\tilde{C}C_{g}C_{t}(m-i)!K^{m-i}\\
       &&( \textrm{ by symmetry of } \overline{q}_{m-i}(t-s,y,z), \textrm{ Proposition
       3.7 and } \int_{\R^d}\overline{p}(s,x,w)\,m(dw)\le D_{2} )\\
       &&+\frac{L^m}{(1-\frac{1}{10}-
       2^{-\frac{1}{2}})^{d+\alpha}}\overline{C}D^2_{2}\frac{\frac{t}{2}}{|x-z|^{d+\alpha}}.\\
\end{eqnarray*}
It is easy to see that there exists $t_{33}>0$ and $t_{33} \leq
\min(t_0,t_2)$ such that when $ 0< t \le t_{33}$, the first three
terms in above inequality $\le
\frac{1}{4}\tilde{C_{1}}m!K^m\frac{t}{|x-z|^{d+\alpha}},\,\,\textrm{for
all } m>0$. We can also have the fourth term in the above
inequality $\le
\frac{1}{8}\tilde{C_{1}}m!K^m\frac{t}{|x-z|^{d+\alpha}},\,\,
\textrm{for all } m>0$, by the choice of $\tilde{C_2}$ and
$\tilde{C_{1}} \ge \tilde{C_{2}}$. Thus in case 2 when $ 0< t \le
t_{33},$
\begin{eqnarray*}
&&\sum_{i=1}^mC^{i}_m\int_{0}^{\frac{t}{2}}\int_{\R^d}\overline{p}(s,x,w)m(dw)\int_{\R^d}\frac{\overline{C}\overline{F}^{i}(w,y)}{|w-y|^{d+\alpha}}\overline{q}_{m-i}(t-s,y,z)\,m(dy)ds\\
        &\le&
        \frac{1}{2}\tilde{C_{1}}m!K^{m}\frac{t}{|x-z|^{d+\alpha}}.
\end{eqnarray*}
Combining case 1 and case 2, when $ 0< t \le \min(t_{23},t_{33}),$
\begin{eqnarray*}
&&\sum_{i=1}^mC^{i}_m\int_{0}^{\frac{t}{2}}\int_{\R^d}\overline{p}(s,x,w)m(dw)\int_{\R^d}\frac{\overline{C}\overline{F}^{i}(w,y)}{|w-y|^{d+\alpha}}\overline{q}_{m-i}(t-s,y,z)\,m(dy)ds\\
        &\le& \frac{1}{2}\tilde{C_{1}}m!K^{m}t^{-\frac{d}{\alpha}}\left(1 \wedge
        \frac{t^{\frac{1}{\alpha}}}{|x-z|}\right)^{d+\alpha}.
\end{eqnarray*}

For the second term in the expression of
$\overline{q}_{m}(t,x,z)$:
$$\sum_{i=1}^mC^{i}_m\int_{\frac{t}{2}}^{t}\int_{\R^d}\overline{p}(s,x,w)m(dw)\int_{\R^d}\frac{\overline{C}\overline{F}^{i}(w,y)}{|w-y|^{d+\alpha}}\overline{q}_{m-i}(t-s,y,z)\,m(dy)ds.$$
Letting $t-s=\tilde{s}$, the second term becomes
$$\sum_{i=1}^mC^{i}_m\int_{0}^{\frac{t}{2}}\int_{\R^d}\overline{p}(t-\tilde{s},x,w)m(dw)\int_{\R^d}\frac{\overline{C}\overline{F}^{i}(w,y)}{|w-y|^{d+\alpha}}\overline{q}_{m-i}(\tilde{s},y,z)\,m(dy)d\tilde{s}.$$
There are two cases.

Case a. When $|x-z| \le t^{\frac{1}{\alpha}}$,
\begin{eqnarray*}
        &&\sum_{i=1}^mC^{i}_m\int_{0}^{\frac{t}{2}}\int_{\R^d}\overline{p}(t-\tilde{s},x,w)m(dw)\int_{\R^d}\frac{\overline{C}\overline{F}^{i}(w,y)}{|w-y|^{d+\alpha}}\overline{q}_{m-i}(\tilde{s},y,z)\,m(dy)d\tilde{s}\\
        &\le& \sum_{i=1}^mC^{i}_m\int_{0}^{\frac{t}{2}}\int_{\R^d}\tilde{C_{1}}(\frac{t}{2})^{-\frac{d}{\alpha}}m(dw)\int_{\R^d}\frac{\overline{C}\overline{F}^{i}(w,y)}{|w-y|^{d+\alpha}}\overline{q}_{m-i}(\tilde{s},y,z)\,m(dy)d\tilde{s}\\
        &\le&
        \sum_{i=1}^mC^{i}_m\tilde{C_{1}}(\frac{t}{2})^{-\frac{d}{\alpha}}L^{i-1}{\overline{M}}^2\overline{C}\int_{0}^{\frac{t}{2}}\int_{\R^d}\overline{q}_{m-i}(\tilde{s},y,z)\,\mu(dy)d\tilde{s}\\
        &\le&
        \sum_{i=1}^mC^{i}_m\tilde{C_{1}}(\frac{t}{2})^{-\frac{d}{\alpha}}L^{i-1}{\overline{M}}^2\overline{C}C_{t}(m-i)!K^{m-i}\\
        && (\textrm{ by symmetry of } \overline{q}_{m-i}(\tilde{s},y,z)
        \textrm{ and Theorem 3.4} )\\
        &=&m!\sum_{i=1}^m(\frac{L^{i-1}K^{m-i}}{i!}){\overline{M}}^2\overline{C}(\frac{1}{2})^{-\frac{d}{\alpha}}C_{t}\tilde{C_{1}}t^{-\frac{d}{\alpha}}.
\end{eqnarray*}
  Since there exists $ \tilde{t}_{23} >0 \textrm{ and } \tilde{t}_{23} \le \min(t_0,t_{13}) $ such that when $0 < t \le
  \tilde{t}_{23}$,
       $$\sum_{i=1}^n(\frac{L^{i-1}K^{n-i}}{i!}){\overline{M}}^2\overline{C}(\frac{1}{2})^{-\frac{d}{\alpha}}C_{t}\le
       \frac{1}{2}K^{n},\,\,\forall n \ge 1,$$
  we have in case a when $0 < t \le \tilde{t}_{23}$,
\begin{eqnarray*}
        &&\sum_{i=1}^mC^{i}_m\int_{0}^{\frac{t}{2}}\int_{\R^d}\overline{p}(t-\tilde{s},x,w)m(dw)\int_{\R^d}\frac{\overline{C}\overline{F}^{i}(w,y)}{|w-y|^{d+\alpha}}\overline{q}_{m-i}(\tilde{s},y,z)\,m(dy)d\tilde{s}\\
        &\le&
        \frac{1}{2}\tilde{C_{1}}m!K^{m}t^{-\frac{d}{\alpha}},
\end{eqnarray*}
i.e.
\begin{eqnarray*}
        &&\sum_{i=1}^mC^{i}_m\int_{\frac{t}{2}}^{t}\int_{\R^d}\overline{p}(s,x,w)m(dw)\int_{\R^d}\frac{\overline{C}\overline{F}^{i}(w,y)}{|w-y|^{d+\alpha}}\overline{q}_{m-i}(t-s,y,z)\,m(dy)ds\\
        &\le&
        \frac{1}{2}\tilde{C_{1}}m!K^{m}t^{-\frac{d}{\alpha}}.
\end{eqnarray*}

Case b. When $|x-z|\ge t^{\frac{1}{\alpha}}$. Let
$\tilde{B_1}=\{y\in\R^d|\,|y-z| \ge \frac{1}{10}|x-z|\},
\tilde{B_2}=\{w\in\R^d|\,|w-x| \ge  2^{-\frac{1}{2}}|x-z|\}$ and
$\tilde{B_3}=\{(w,y)\in \R^d \times \R^d|\,|y-z| <
\frac{1}{10}|x-z|,\, |w-x| <  2^{-\frac{1}{2}}|x-z|\} $. On
$\tilde{B_3}$, we have $|w-y| \ge (1-\frac{1}{10}-
2^{-\frac{1}{2}} )|x-z|$.
\begin{eqnarray*}
       &&\sum_{i=1}^mC^{i}_m\int_{0}^{\frac{t}{2}}\int_{\R^d}\overline{p}(t-\tilde{s},x,w)m(dw)\int_{\R^d}\frac{\overline{C}\overline{F}^{i}(w,y)}{|w-y|^{d+\alpha}}\overline{q}_{m-i}(\tilde{s},y,z)\,m(dy)d\tilde{s}\\
       &\le&\sum_{i=1}^mC^{i}_m\int_{0}^{\frac{t}{2}}\int_{\R^d}\overline{p}(t-\tilde{s},x,w)m(dw)\int_{\R^d}\frac{\overline{C}\overline{F}^{i}(w,y)}{|w-y|^{d+\alpha}}\overline{q}_{m-i}(\tilde{s},y,z)1_{\tilde{B_{1}}}(y)\,m(dy)d\tilde{s}\\
       &&+\sum_{i=1}^mC^{i}_m\int_{0}^{\frac{t}{2}}\int_{\R^d}\overline{p}(t-\tilde{s},x,w)m(dw)\int_{\R^d}\frac{\overline{C}\overline{F}^{i}(w,y)}{|w-y|^{d+\alpha}}\overline{q}_{m-i}(\tilde{s},y,z)1_{\tilde{B_{2}}}(w)  \,m(dy)d\tilde{s}\\
       &&+\sum_{i=1}^mC^{i}_m\int_{0}^{\frac{t}{2}}\int_{\R^d}\overline{p}(t-\tilde{s},x,w)m(dw)\int_{\R^d}\frac{\overline{C}\overline{F}^{i}(w,y)}{|w-y|^{d+\alpha}}\overline{q}_{m-i}(\tilde{s},y,z)1_{\tilde{B_{3}}}(w,y) \,m(dy)d\tilde{s}\\
       &\le& \sum_{i=1}^mC^{i}_m\int_{0}^{\frac{t}{2}}\int_{\R^d}\overline{p}(t-\tilde{s},x,w)m(dw)\int_{\R^d}\frac{\overline{C}\overline{F}(w,y)}{|w-y|^{d+\alpha}} \overline{M}L^{i-1}\tilde{C_{1}}(m-i)!K^{m-i}\frac{10^{d+\alpha}\tilde{s}}{|x-z|^{d+\alpha}}\,dyd\tilde{s}\\
       &&+\sum_{i=1}^mC^{i}_m\int_{0}^{\frac{t}{2}}\int_{\R^d}\,m(dw)\tilde{C_{1}}\overline{M}2^{\frac{1}{2}(d+\alpha)}\frac{(t-\tilde{s})}{{|x-z|}^{d+\alpha}}\int_{\R^d}\frac{\overline{C}\overline{F}(w,y)}{|w-y|^{d+\alpha}}\overline{q}_{m-i}(\tilde{s},y,z)\,dydsL^{i-1}\\
       &&+\sum_{i=1}^{m-1}C^{i}_m\frac{L^i}{(1-\frac{1}{10}- 2^{-\frac{1}{2}})^{d+\alpha}}\overline{C}\int_{0}^{\frac{t}{2}}\int_{\R^d}\overline{p}(t-\tilde{s},x,w)m(dw)\int_{\R^d}\frac{1}{|x-z|^{d+\alpha}}\overline{q}_{m-i}(\tilde{s},y,z)\,\\
       &&\cdot m(dy)d\tilde{s}+\frac{L^m}{(1-\frac{1}{10}-2^{-\frac{1}{2}})^{d+\alpha}}\overline{C}D^2_{2}\frac{\frac{t}{2}}{|x-z|^{d+\alpha}}\\
\end{eqnarray*}
\begin{eqnarray*}
       &\le& \sum_{i=1}^mC^{i}_mC_{t}
       {\overline{M}}^2\overline{C}L^{i-1}\tilde{C_{1}}(m-i)!K^{m-i}10^{d+\alpha}\frac{t}{|x-z|^{d+\alpha}}\\
       &&+\sum_{i=1}^mC^{i}_m\tilde{C_{1}}{\overline{M}}^2\overline{C}2^{\frac{1}{2}(d+\alpha)}\frac{t}{{|x-z|}^{d+\alpha}}C_{t}(m-i)!K^{m-i}L^{i-1}\\
       &&( \textrm{ by symmetry of } \overline{q}_{m-i}(\tilde{s},y,z)\textrm{ and Theorem 3.4 } )\\
       &&+\sum_{i=1}^{m-1}C^{i}_m\frac{L^i}{(1-\frac{1}{10}- 2^{-\frac{1}{2}})^{d+\alpha}}\overline{C}D\frac{t}{|x-z|^{d+\alpha}}\tilde{C}C_{g}C_{t}(m-i)!K^{m-i}\\
       &&( \textrm{ by symmetry of } \overline{q}_{m-i}(\tilde{s},y,z), \textrm{ Proposition
       3.7 and } \int_{\R^d}\overline{p}(t-\tilde{s},x,w)\,m(dw)\le D_2)\\
       &&+\frac{L^m}{(1-\frac{1}{10}-
       2^{-\frac{1}{2}})^{d+\alpha}}\overline{C}D^2_{2}\frac{\frac{t}{2}}{|x-z|^{d+\alpha}}.\\
\end{eqnarray*}
It is easy to see that there exists $\tilde{t}_{33}>0$ and
$\tilde{t}_{33}\leq\min(t_0,t_2)$ such that when $ 0< t \le
\tilde{t_{33}}$, the first three terms in above inequality $\le
\frac{1}{4}\tilde{C_{1}}m!K^m\frac{t}{|x-z|^{d+\alpha}},\,\,\textrm{for
any } m>0$. We can also have the fourth term in the above
inequality $\le
\frac{1}{8}\tilde{C_{1}}m!K^m\frac{t}{|x-z|^{d+\alpha}},\,\,\textrm
{for any } m>0$, by the choice of $\tilde{C_2}$ and $\tilde{C_{1}}
\ge \tilde{C_{2}}$. Thus in case b when $ 0< t \le
\tilde{t}_{33},$
\begin{eqnarray*}
        &&\sum_{i=1}^mC^{i}_m\int_{0}^{\frac{t}{2}}\int_{\R^d}\overline{p}(t-\tilde{s},x,w)m(dw)\int_{\R^d}\frac{\overline{C}\overline{F}^{i}(w,y)}{|w-y|^{d+\alpha}}\overline{q}_{m-i}(\tilde{s},y,z)\,m(dy)d\tilde{s}\\
        &\le&
        \frac{1}{2}\tilde{C_{1}}m!K^{m}\frac{t}{|x-z|^{d+\alpha}},
\end{eqnarray*}
i.e.
\begin{eqnarray*}
&&\sum_{i=1}^mC^{i}_m\int_{\frac{t}{2}}^{t}\int_{\R^d}\overline{p}(s,x,w)m(dw)\int_{\R^d}\frac{\overline{C}\overline{F}^{i}(w,y)}{|w-y|^{d+\alpha}}\overline{q}_{m-i}(t-s,y,z)\,m(dy)ds\\
        &\le&
        \frac{1}{2}\tilde{C_{1}}m!K^{m}\frac{t}{|x-z|^{d+\alpha}}.
\end{eqnarray*}
Combining case a and case b, when $ 0< t \le
\min(\tilde{t}_{23},\tilde{t}_{33}),$
\begin{eqnarray*}
&&\sum_{i=1}^mC^{i}_m\int_{\frac{t}{2}}^{t}\int_{\R^d}\overline{p}(s,x,w)m(dw)\int_{\R^d}\frac{\overline{C}\overline{F}^{i}(w,y)}{|w-y|^{d+\alpha}}\overline{q}_{m-i}(t-s,y,z)\,m(dy)ds\\
        &\le& \frac{1}{2}\tilde{C_{1}}m!K^{m}t^{-\frac{d}{\alpha}}\left(1 \wedge
        \frac{t^{\frac{1}{\alpha}}}{|x-z|}\right)^{d+\alpha}.
\end{eqnarray*}
Therefore when $0 < t <
t_{3}=\min(t_{23},t_{33},\tilde{t}_{23},\tilde{t}_{33})$,
$$ \overline{q}_{m}(t,x,z) \le \tilde{C_{1}}m!K^{m}t^{-\frac{d}{\alpha}}\left(1 \wedge
\frac{t^{\frac{1}{\alpha}}}{|x-z|}\right)^{d+\alpha},$$ i.e. the
statement holds for $n=m$. \qed

By Theorem 3.8, we have for $0<t\leq t_{3}$,
$$ \sum_{n=0}^{\infty}\frac{\overline{q}_n(t,x,z)}{n!}\leq \sum_{n=0}^{\infty}\tilde{C_{1}}K^nt^{-\frac{d}{\alpha}}\left(1 \wedge \frac{t^{\frac{1}{\alpha}}}{|x-z|}\right)^{d+\alpha}  =\tilde{C_{1}}\frac{1}{1-K}t^{-\frac{d}{\alpha}}\left(1 \wedge \frac{t^{\frac{1}{\alpha}}}{|x-z|}\right)^{d+\alpha}. $$
Since $|q_n(t,x,z)|\le \overline{q}_n(t,x,z)$ and
$q(t,x,z)=\sum_{n=0}^{\infty}(-1)^n\frac{{q}_n(t,x,z)}{n!}$,
$$q(t,x,z)\le
\tilde{C_{1}}\frac{1}{1-K}t^{-\frac{d}{\alpha}}\left(1 \wedge
\frac{t^{\frac{1}{\alpha}}}{|x-z|}\right)^{d+\alpha}.$$

Applying (ii) of Proposition 3.5, we have
$$q(t,x,z)\le {C}_{5}e^{{C}_{6}t}t^{-\frac{d}{\alpha}}\left(1 \wedge \frac{t^{\frac{1}{\alpha}}}{|x-z|}\right)^{d+\alpha},\,\,\forall (t,x,y)\in (0,\infty) \times \R^d \times \R^d,$$
where $C_{5}$ and $C_{6}$ are positive constants. Thus we have
established the lower and upper estimates of $q(t,x,y)$ as
follows,
\begin{thm}
There exist positive constants $C_{3}$,$C_{4}$,$C_{5}$ and $C_{6}$
such that
\begin{equation}
C_{3}e^{-C_{4}t}t^{-\frac{d}{\alpha}}\left(1 \wedge
\frac{t^{\frac{1}{\alpha}}}{|x-z|}\right)^{d+\alpha} \leq q(t,x,z)
\leq  C_{5}e^{C_{6}t}t^{-\frac{d}{\alpha}}\left(1 \wedge
\frac{t^{\frac{1}{\alpha}}}{|x-z|}\right)^{d+\alpha}
\end{equation} for all $(t,x,z)\in (0,\infty) \times \R^d \times
\R^d $.

\end{thm}

From Theorem 3.9 and (ii) of Proposition 3.5, it is easy to obtain
the following property for $q(t,x,z)$,
\begin{prop}
 $q(t,x,z)$ is joint continuous on
$(0,\infty)\times\R^d\times\R^d.$
\end{prop}

\bigskip
\bigskip
{\bf Acknowledgement}: I am very grateful to my advisor Professor
Renming Song for his encouragement and many suggestions.

\vspace{.5in}
\begin{singlespace}

\end{singlespace}

\end{doublespace}

\end{document}